\documentclass[11pt,reqno]{amsart}
\usepackage{amsmath}
\usepackage{amsthm}
\usepackage{graphicx}
\usepackage{hyperref}
\usepackage{mathrsfs}
\usepackage{cite}
\usepackage{dsfont}
\usepackage{color}
\usepackage{bm}
\usepackage{amsfonts}
\usepackage{amssymb}
\usepackage{enumerate}
\usepackage{soul}
\usepackage[margin=1.1in]{geometry}
\usepackage[normalem]{ulem}

\makeatletter
\def\@settitle{\begin{center}%
		\baselineskip14\p@\relax
		\bfseries
		\uppercasenonmath\@title
		\@title
		\ifx\@subtitle\@empty\else
		\\[1ex]\uppercasenonmath\@subtitle
		\footnotesize\mdseries\@subtitle
		\fi
	\end{center}%
}
\def\subtitle#1{\gdef\@subtitle{#1}}
\def\@subtitle{}
\makeatother

\numberwithin{equation}{section}

\newtheorem{theorem}{Theorem}[section]
\newtheorem{lemma}[theorem]{Lemma}
\newtheorem{proposition}[theorem]{Proposition}
\newtheorem{corollary}[theorem]{Corollary}

\theoremstyle{definition}

\newtheorem{definition}[theorem]{Definition}
\newtheorem{remark}[theorem]{Remark}

\definecolor{darkgreen}{rgb}{0,0.5,0}
\definecolor{darkblue}{rgb}{0,0,0.7}
\definecolor{darkred}{rgb}{0.9,0.1,0.1}
\definecolor{lightblue}{rgb}{0,0.51,1}

\newcommand{\R}{{\mathbb R}}
\newcommand{\Rd}{{\mathbb{R}^d}}
\newcommand{\N}{{\mathbb N}}
\newcommand{\E}{{\mathds E}}
\newcommand{\Pb}{{\mathds P}}
\newcommand{\F}{{\mathcal F}}
\newcommand{\Q}{{\mathcal Q}}

\newcommand{\dx}{\partial_x}
\newcommand{\dxx}{{\partial^2_{xx}}}

\newcommand{\1}{\mathds{1}}

\newcommand{\ds}{\displaystyle}

\DeclareMathOperator{\tr}{Tr}

\newcommand{\xit}{X^{i,N}_t}
\newcommand{\xjt}{X^{j,N}_t}
\newcommand{\xis}{X^{i,N}_s}
\newcommand{\xjs}{X^{j,N}_s}

\newcommand{\Bit}{B^i_t}
\newcommand{\Bot}{B^0_t}
\newcommand{\Bis}{B^i_s}

\newcommand{\eit}{E^{i,N}_t}
\newcommand{\zit}{Z^{i,N}_t}
\newcommand{\ejt}{E^{j,N}_t}
\newcommand{\eis}{E^{i,N}_s}
\newcommand{\ejs}{E^{j,N}_s}
\newcommand{\eio}{E^{i,N}_0}

\newcommand{\eism}{E^{i,N}_{s^-}}
\newcommand{\ejsm}{E^{j,N}_{s^-}}

\newcommand{\mun}{\mu^N_t}

\newcommand{\mut}{\mu_t}

\title[Conditional propagation of chaos in a spatial stochastic epidemic model]{Conditional propagation of chaos in a spatial stochastic epidemic model with common noise}
\author{Yen V. Vuong, Maxime Hauray, Étienne Pardoux}

\begin{document}
\begin{abstract}
We study a stochastic spatial epidemic model where the $N$ individuals carry two features: a position and an infection state, interact and move in $\R^d$. In this Markovian model, the evolution of the infection states are described with the help of the Poisson Point Processes , whereas the displacement of the individuals are driven by mean field advection, a (state dependence) diffusion and also a common noise, so that the spatial dynamic is a random process. We prove that when the number $N$ of individual goes to infinity, the conditional propagation of chaos holds : conditionnally to the common noise, the individuals are asymptotically independent and the stochastic dynamic converges to a "random" nonlinear McKean-Vlasov process. As a consequence, the associated empirical measure converges to a measure, which is solution of a stochastic mean-field PDE driven by the common noise.
\end{abstract}
	
\maketitle
	
	
	%\tableofcontents
%%%%%%%%%%%%%%%%%%%%%%%%%%%%%%%%%%%%%%
%
%
%%%%%%%%%%%%%%%%%%%%%%%%%%%%%%%%%%%%%%
\section{Introduction}

Epidemic models have been studied for a long while, in both deterministic and stochastic settings. 
%By law of large number,  deterministic model can be obtained as the limit of a stochastic epidemic model, see \cite{Britton-Pardoux}. 
In this paper, we study a spatial model  based on the famous SIR model, the letter S, I and R standing for the different states of an individual which can pass from the compartment of "Susceptible" to the "Infected" one and eventually to the compartment of "Recovered" when the individual recovers from the illness.
In our spatial model, an individual will be characterized by:
\begin{itemize}
\item Its state $E \in  \{S,I,R\} = \{ 0,1,2\}  $, since  we will identify $S$ with $0$, $I$ with $1$ and $R$ with $2$ in order to simplify the mathematical description,
\item Its position, a continuous variable $X \in \R^d$. Typically, $d=1,2$ for the propagation of epidemics\footnote{but the case $d=3$ is also considered in quite similar models  in chemistry for the microscopical description of chemical reactions, see~\cite{Lim-Lu-Nolen}.}.
\end{itemize} 

So the individual phase space is $\Pi := \R^d \times \{0,1,2\}$.  We consider a community of $N$ individuals, denoting by  $\zit = \bigl(\xit,\eit)$, $i=1,\ldots,N$, the position and state of the $i^{\text{th}}$ individual at time $t >0$. The full vector $(\zit)_{i \le N}$ belongs to the full phase space $\Pi^N = \R^{dN} \times \{0,1,2\}^N$.

The introduction of spatial variables will complicate the standard homogeneous SIR model in two directions: by using an infection rate that depends on the position of the individuals, and by taking into account the individual displacements.

\medskip
\paragraph{\bf An infection rate depending on the position}

% can use  infection rate of an individual from the others that depends on their positions (for instance their relative distance). 

It is quite clear that in realistic situation, an infected individual will infect a close neighbour with a higher rate than an individual living far away. While these different behaviours are averaged in a homogeneous SIR model, in spatial models we use an infection rate depending on the positions. The infection rate between locations $x$ and  $y \in \R^d$ will be given by a function $K : \R^d \times \R^d \to  \R_+$, that we will assume to be bounded and Lipschitz. Averaging over all the infected individuals, the susceptible individual $i$ becomes infected (in other words its state jumps from $0$ to $1$) at time $t$ at the rate
\begin{equation} \label{rate_infection}
\frac1{N}\sum_{j=1}^N K(\xit,\xjt)\1_{\{\eit=S\}}\1_{\{\ejt=I\}}.
\end{equation}

What could be the form of a realistic function $K$?
For several sedentary species, it could be reasonnable to think that the individuals could only infect their close neighbours, and use a compactly supported and constant kernel, of the form
\[
K(x,y) = \frac{\beta}{c_d r_0^d} \1_{|x-y| \le r_0}
\]
where $r_0$ is the size of the "territory" of the individuals, and $c_d$ is the volume of the unit ball, and $\beta$ is the infection rate.

For the propagation of an epidemic among a human population, it is possible to incorporate in $K$ effects due to the human mobility, which has become a subject of intensive research in the past decades. Most popular models for the human mobility are the so-called gravity model~\cite{erlander} which state that the total number of travel between two cities is inversely proportional to the inverse of the square of the distance between them\footnote{Some versions allow to fit the power of the distance rather than using the power $2$, or even to replace it by a more general function of the distance.}, and the radiation model~\cite{simini}, which state that the same number depends also on the distribution of population between the two locations.
The gravity model naturally leads in our case in the use of a function $K$ of the form 
\[
K(x,y) = \frac{c}{d^2 + |x-y|^2},
\]
with $c,d \in \R^{+*}$ two constants to be fitted with the observations. According to the study~\cite{RBBSB}, that choice is in a good accordance with the observed data on the spread of the Covid epidemic in France in 2020.
We refer to~\cite{BGVK} for the description of other mobility functions and the study of their relevance in (discrete) spatial epidemic model. To be complete, we also mention a recent study~\cite{Schlapfer} that deduces from the observations a quite general law in urban mobility and that could lead to a general kernel, not only depending on the relative distance between $x$ and $y$, but also on the atractivity of the different locations. The choice of an adapted kernel is clearly a relevant question, which is still a subject under investigation, and here we will consider a general kernel, depending on the two positions.

%\mh{$\beta$ was useless as it could be integrate in $K$. So I erase it here and in the sequel.}

%\vvy{But as usual, the quantity $R_0:=\frac{\beta}{\gamma}$ always plays an important role in the epidemic models. It is called the reproduction number. For example, in the model without spatial structure, this quantity controls the epidemic, which could blow up or decay. I am not sure in general whether we need it or not...}

%\mh{ Your are right in a homogeneous model. But here in that spatial model, the $R_0$ is not defined by $\beta/ \gamma$, at least with the $\beta$ that was define before. Now the $\beta$ is local, something like $\beta(x) = \int K(x,y) \mu(dy)$. And the global version of $\beta$, if it exists is something like $\iint K(x,y) \mu(dx,0) \mu(dy,0) /  \mu(0)$.  We could write a remark about that.}

\medskip
\paragraph{\bf Recovery rate}
 The infectious individuals recover (in other words their state jumps from $1$ to $2$) at rate $\gamma >0$ independently of everything else (of the other individuals, the number of infected and the respective positions of the individuals,...) and once an individual recovers, it can not become infected anymore.

\medskip
\paragraph{\bf Movements of individuals}

In order to study the propagation of epidemics among non sedentary species, we will take into account some kind of displacements of the individuals. In our model, each individual  moves in $\Rd$ according to 
\begin{itemize}
     \item A mean field interaction (with all the  others individuals),
     \item Independent diffusions,
     \item A random drift (or diffusion), common to all the individuals.
\end{itemize}
So the evolution of the position $\xit$ of the $i^{\text{th}}$ individual satisfies the following equation
\begin{equation} \label{movement}
    d\xit=\displaystyle\frac{1}{N}\sum_{j=1}^{N}V \Bigl(\xit,\eit,\xjt,\ejt \Bigr)dt+\sigma\Bigl(\xit,\eit\Bigr)d\Bit+\sigma_0\bigl(\xit,\eit\bigr)d\Bot
\end{equation}
where the interaction kernel $V$, the diffusion strength $\sigma$ and the strengh $\sigma_0$ of the random drift-diffusion all depend on the individual's state and position, and are bounded Lipschitz continuous with respect to the position variables. Of course, this equation has a meaning on a probability space endowed with the requested Brownian motions (and Poisson point processes for the jump-infectious part of the dynamic). The Lipschitz hypothesis will be very usefull to build a correct theory of exsitence, uniqueness to that system, and also for our results concerning the large population limit (i.e. when $N$ goes to infinity).

\medskip
Before going on, let us discuss a bit the possible application of that kind of dynamics, and the particular choice of diffusion, mean-field and common drift-diffusion that could be relevant.

Mean-field interactions are commonly chosen in order to describe swarming of flocking behaviour of some species, see for instance~\cite{Carrillo}. Such phenomena occur on time scales shorter than the one of epidemic propagation. We choose to use it in order to illustrate that our model is compatible with the introduction of a kind of collective interaction.

The common drift term is here written in a quite simple form. It could probably be generalized to a more general random drift without diffusion which is probably the most interesting case) like in~\cite{Coghi-Flandoli} of the form
\[
\sum_{k=1}^\infty \sigma_k(\xit,\eit) dB_t^{0,k},
\]
with divergence-free $\sigma_k$, which is a quite general way to describe homogeneous random drift\footnote{It is shown in~\cite{Coghi-Flandoli} that under the divergence free hypothesis, the Ito and Stratonovich formulation are equivalent, so that we have the right to call it random drift here.}.

The individual diffusion term is quite common to model the displacement of individuals.  Modeling it with the help of a standard diffusion is a first step, that is probably not very realistic, since individuals usually do not have a Markovian behaviour. One more realistic model among human population is for instance the EPR model (Exploration and Preferential Return) introduced in~\cite{song}.

We stress out the fact that the diffusion strength depends on the infectious state of the individuals is quite interesting for the application. We could think of the propagation of a rabies epidemic, where the disease affects the displacement behaviour of the individuals.
But also in the case of epidemics in human population, as the individuals are exchangeable, we may study the case where the susceptible population is at rest, while the infectious individuals follow a diffusion in order to model the dispersion of the epidemic. In that case, the classical diffusion could be an interesting choice, even if some more realistic models should probably include Lévy flights and also be non Markovian, as suggested by the study of banknote diffusion in a population~\cite{Brockmann}.

%\cite{simini} seems useless. To check.
%Random behaviour in a random environment, with the EPR model as suggested in~\cite{Schlapfer}.
%	The individual noises $(\Bit)_{t\geq 0},\;i=0,\dotsc,N$ are independent $\Rd$-valued Brownian motions, and the common noise $\Bot$ is also a $\Rd$-valued Brownian motion independent of the previous ones.

%\mh{Il faudra une discussion sur ce que veut dire une fonction $V: \Pi^2 \to \R^d$ Lipschitz. C'est un peu bizarre de parler de Lipschitz par rapport à $E$ qui est une variable discrète. Après des calculs, je vois que la constante de Lipschitz de V par rapport à $(x,e,y,f)$ vérifie :
%\[
%L_V  \le 
%\max_{i,j = 0,1,2} \bigl( \| \partial_x V_{i,j}\|_\infty +
%\| \partial_y V_{i,j}\|_\infty \bigr) + 
%\max_{i,i',j,j' = 0,1,2} \| V_{i,j} - V_{i',j}\|_\infty
%\]
%avec $V_{i,j}(x,y) = V(x,i,y,j)$.
%}
\subsection{The model}
In view of the above settings, a description of the epidemiological dynamic with the help of Poisson Point processes is suitable. If we choose a probability space equipped with $N$ independent Poisson point processes $(P^i)_{i=1,\dots,N}$ and $N+1$ Brownian motions, the position and state of the individuals will evolve according to the following dynamics :
\begin{equation}\label{system}
    \begin{cases}
    d\xit & \hspace{-.4cm} = \displaystyle\frac{1}{N}\sum_{j=1}^{N}V \bigl(\xit,\eit,\xjt,\ejt \bigr)dt+\sigma\bigl(\xit,\eit\bigr)d\Bit+\sigma_0\bigl(\xit,\eit\bigr)d\Bot,
    \\
    \eit & \hspace{-.4cm} = \eio + P^i \biggl( \displaystyle\int_0^t   \Bigl\{ \frac1{N} \sum_{j \neq i} K\bigl( \xis,\xjs  \bigr) \1_{(\eism,\ejsm)=(0,1)} + \gamma \1_{\eism=1} \Bigr\} \,ds \biggr).
    \end{cases}
\end{equation}

This model is inspired by the previous work~\cite{Emakoua} where the case without common random drift and mean-field interaction is treated, and results about the Law of Large Number and the Central Limit Theorem are obtained. Here we will generalize the Law of Large Number result.

%This work is inspired by \cite{Emakoua}, where the propagation of an epidemic is also studied in a spatial configuration but for a simpler dynamic: the individuals follow  independent Brownian motions whereas the epidemic dynamic is the same. 
%Our spatial dynamic is more complex: the displacement of individuals involves  a mean field terms and the diffusion depends on their own position and states. 

Let us also mention some previous works in other related contexts where a common noise  affects  all individuals :  interacting particle systems  with a common random drift \cite{Coghi-Flandoli},  mean field games \cite{Carmona-Delarue}, \cite{Carmona-Delarue-Lacker}.

An important concept in studing the mean field limit of the interacting pacticle system is the propagation of chaos \cite{Sznitman, Hauray-Mischler}, which roughly states that when the number of pacticles tends to infinity, the chaotic character of the initial distribution propagates through the system at any time. In the present model, the individuals always keep track of the common noise $\Bot$ so it is impossible to expect the asymptotically independence at the limit. Nevertheless, it is reasonable to expect that the individuals become asymptotically independent conditionally on the information generated by the common noise. The objective of this paper is to establish the conditional propagation of chaos of the $N$ individuals process towards a random nonlinear (or McKean-Vlasov) process defined later. This also implies the convergence of the empirical measure process
\[
t  \mapsto \mun=\frac{1}{N}\sum_{i=1}^{N}\delta_{(\xit,\eit)}
\]
which is a random process with value in  $\mathcal P(\Pi):=\mathcal P \bigl(\Rd\times\{0,1,2\}\bigr)$. The proof will follow a standard coupling method introduced by McKean in \cite{McKean}.

\subsection{The existing literature and the novelty of our result.}

The idea of our paper was inspired originally mostly by the two papers~\cite{Emakoua} and~\cite{Coghi-Flandoli}. Here we detail the novelties and improvements in our work with respect to these two references.

\medskip
Our first inspiring paper~\cite{Emakoua} was written  by A. Emakoua, S. Bowong and the third author. They studied a spatial epidemic model which is similar to ours but a bit simpler, without 
\begin{itemize}
    \item the common stochastic case (with our notation $\sigma_0=0$) ,
    \item interaction among the individuals (with our notation $V=0$).
\end{itemize}
In their model, the positions of the individuals are independent Brownian motions. In that quite simpler setting, their main result is a Central Limit Theorem (CLT) in the limit of a large number of individuals. The  proof of the CLT result is delicate and technical, and relies of course on a Law of Large Number (LLN) result, much simpler to prove. In our present article we focus on the  LLN result, and improve it in two directions:
\begin{itemize}
    \item We prove it  in a more general (and realistic) situation, 
    \item We obtain \emph{quantitative} estimates: in Wasserstein distance rather than a qualitative result stating the weak convergence. 
\end{itemize}
Remark that the stronger CLT result of~\cite{Emakoua} also gives as a byproduct of the order of convergence in the LLN, but not for the Wasserstein distance as we obtain. We also mention that the proof of a CLT in our setting will be the aim of a future work of the first author.

\medskip
Our second inspiring article was written by M. Coghi and F. Flandoli \cite{Coghi-Flandoli}. In that work, the authors prove a conditionnal propagation of chaos result for an interacting particle system where the particles (having only a position parameter and no distinct types) are subject to individual and common noises. With our formalism, there are two main differences:
\begin{itemize}
    \item  There is no internal state $E$: so the dynamical equation in their model is very similar to the first equation in~\eqref{system}, used with functions $V,\sigma,\sigma_0$ independent of $E$. \item They used   a quite general (smooth enough) cylindrical BM with quite general (smooth) correlation function. That noise could be written as a infinite serie of the form $dW(s,x) = \sum_{k} \sigma^k(x) dB^{0,k}_s$, where $(B^{0,k}_s)_{s \ge 0}$ are standard Brownian motions and $(\sigma^k)_{k \ge 0}$ a family of smooth enough functions. Their common noises are ``richer'' and more realistic than ours since we use only \emph{one} standard BM to define ours. But it seems reasonable to expect that our result extend to the more realistic situation studied in~\cite{Coghi-Flandoli}.
\end{itemize}

If we restrict our result to the case where the movement does not depend upon the internal state ($V$,$\sigma$,$\sigma_0$ independent of $E$) so that we can study the dynamics of the position only, 
and if we use only our simpler type of common noise in~\cite{Coghi-Flandoli}, then 
both results are very similar. We both provide the same  quantitative estimate of conditionnal propagation of chaos in Wasserstein distance of order one. The only difference is in the constants appearing in the estimates. But even if the results are similar, there are still two differences that we would like to  emphasize. 

The first difference is that the use of Wasserstein distance of order one is motivated in~\cite{Coghi-Flandoli} by the willingness to use the standard tool in the field, in order to compare with previous results. But in their setting, the use of the Wasserstein distance of order two, which is defined with a quadratic cost and thus behave much better than the order one distance with respect to the diffusion, would have simplified the proof. In our model, the use of the order one Wasserstein distance is mandatory in order to treat the Poissonian part of the epidemiological dynamics. See the Remark \ref{rmk_using_W1} for a more detailled explanation of that point.

Our proof contains a simplification that could be used also in the work by M. Coghi and F. Flandoli. In order to control the evolution of the Wasserstein distance of order one under diffusion, they used a short time estimate that allows them to iterate in order to get the desired etimate on arbitrary finite time, the iteration step being a bit technical. We perform a more careful analysis of the evolution of the same quantity and get directly an estimate valid for an arbitrary finite time, without any iteration step. See Remark~\ref{rmk_a_direct_estimate} in Section~\ref{sec:proof_NLSDE} for a more precise description of the improvment.

\medskip

\paragraph{\bf Organisation of the paper.}  

We present the model state our main results in Section~\ref{sec:main_results}. The next sections are devoted to the proof of the main results. Section~\ref{sec:proof_NLSDE} is devoted to the proof of Theorem~\ref{exist_unique_thm} and to a result on the propagation of moments. Section~\ref{sec:nlpde} is devoted to the proof of Theorem~\ref{prop:NLSDE_to_SPDE} and some discussions on the uniqueness of solutions to the mean field limit equation. Section~\ref{sec:LLN} contains the proof of Theorem~\ref{Prop_control_couplings} and Corollary~\ref{rate_of_convergence}.

\section{Preliminaries and main results}\label{sec:main_results}
\subsection{Notations}

\medskip
\paragraph{\bf Basic notations}
\begin{itemize}
 \item Variables denoted by $x$ or $y$ (also with capital letters) denotes a spatial position in $\R^d$,
 \item Variables denoted by $e$ (resp. $E$) or $f$ (resp. $F$) denotes an epidemiological status in $\{0,1,2\}$,
 \item Variables denoted by a $z$ (resp. $Z$) denotes a pair $(x,e)$ of a position and an epidemiological status.
\end{itemize}

\medskip
\paragraph{\bf The underlying space}
Throughout this paper, we keep assuming that the independent noises $(B^i)_{i=1,\dots,N}$ are constructed on a probability space $(\widetilde{\Omega},\widetilde{\mathcal{F}},\widetilde{\Pb})$ and the common noise $B^0$ is constructed on another space $(\Omega^0,\mathcal{F}^0,\Pb^0)$. We now define the product structure $(\Omega,\F,(\F_t)_{t\geq 0},\Pb)$, where $\Omega=\widetilde{\Omega}\times\Omega^0$, $(\F,\Pb)$ is the completion of $(\widetilde{\mathcal{F}}\otimes\mathcal{F}^0,\widetilde{\Pb}\otimes\Pb^0)$, and $(\F_t)_{t\geq 0}$ is the complete right continuous version of $(\widetilde{\mathcal{F}}_t\otimes\mathcal{F}^0_t)_{t\geq 0}$. With this construction, if we consider $Z$ a $\Pi$-valued random variable on $(\Omega,\F,\Pb)$ then the conditional law of $Z$ given $\F^0$, denoted by $\mathcal L (Z | \mathcal F^0)$, is a $\mathcal{P}(\Pi)$-valued random variable defined on $(\Omega^0,\mathcal{F}^0,\Pb^0)$.

\medskip
\paragraph{\bf Simplification of the system~\eqref{system}}
With the convenient notation
\[
V^{i,N}_s :=  \frac1N \sum_{j \neq i} V( \xis,\eis,\xjs,\ejs), \qquad
K^{i,N}_s :=  \frac1{N} \sum_{j \neq i} K( \xis,\xjs) \1_{\ejsm=1} 
,
\]
the system~\eqref{system} becomes
\begin{equation}\label{short-system}
    \left\{ \begin{array}{rl}
    %\begin{cases}
    \ds  d\xit & \ds =  \ds V^{i,N}_t dt+\sigma\bigl(\xit,\eit\bigr) d\Bit
    +\sigma_0\bigl(\xit,\eit\bigr)d\Bot,
    \\
    \ds \eit & \ds = \eio + P^i \biggl( \displaystyle\int_0^t   \Bigl\{ K^{i,N}_s \1_{\eism=0} + \gamma \1_{\eism=1} \Bigr\} \,ds \biggr).
    %\end{cases}
    \end{array} \right.
\end{equation}

\medskip
\paragraph{\bf Poisson Random Measure rather than Poisson Point Process}
The counting processes describing the evolution of the type of each individual used in~\eqref{short-system} can be rewritten with the help of independent Poisson random measures  $(\Q^i)_{i=1,\dots,N}$ on $\R_+ \times \R_+$ with mean measure the Lebesgue measure $dt \times du$:
\[
\eit= \eio + 
\int_{[0,t]\times\R_+}\1_{\big\{ u \le K^{i,N}_{s} \1_0(\eism) + \gamma \1_1(\eism)\big\} } \Q^i(ds,du),
\]
In fact the  integral quantity in the r.h.s. has the same law (trajectorially) as the stochastic term in the r.h.s in the second line of~\eqref{short-system}. So the system of equations~\eqref{short-system} can be rewritten as follows
\begin{align}\label{short-system2}
   \left\{
    \begin{array}{rl}
        \ds d\xit &= \ds V^{i,N}_s dt+\sigma\bigl(\xit,\eit\bigr) d\Bit
    +\sigma_0\bigl(\xit,\eit\bigr)d\Bot,  \\
    \ds \eit &= \ds  \eio + \displaystyle\int_{[0,t]\times\R_+}\1_{\big\{ u \le K^{i,N}_{s} \1_0(\eism) + \gamma \1_1(\eism)\big\} } \Q^i(ds,du).
    \end{array} \right.
\end{align}

Even if these random measures $\Q^i$ may seem less simple to handle than the PPPs $P^i$, it will be of a great help in the sequel when we couple processes. Thanks to these PRMs, we will use a true parallel coupling, in the sense that we will construct coupled processes that jump very often exactly at the same time. It is not possible to construct the same coupling by using PPPs. In the later case, our coupled processes will jump at close times but not exactly at the same times.
This has an important consequence: we will be able to obtain estimates of sup norms without allowing for small change in time as it is usual in the Skorokhod topology. This simplifies a lot the proof and allow us to give quantitative convergence results.

\begin{proposition} \label{prop:exis_uniq}
Strong existence and weak uniqueness holds for both system~\eqref{system} and system~\eqref{short-system2}. More precisely, given a probability space endowed with the appropriate Brownian motions and Poisson Point Process in the case of system~\eqref{system} (or Poisson Random Measures in the cases of system~\eqref{short-system2}) , there exists a unique solution to the system~\eqref{system} or~\eqref{short-system2}.

Moreover, solutions to~\eqref{system} or~\eqref{short-system2} constructed on possibly different probability spaces have the same law.
\end{proposition}

The proof of Proposition \ref{prop:exis_uniq} relies on standard argument since all the coefficients appearing in the equation are smooth enough. We refer for instance to~\cite[Theorem 9.1]{Ikeda-Watanabe}.

% \medskip
% \mh{
% \paragraph{\bf Distances used in that article}

% For $z=(x,v)$ and $z'=(x',v') \in \Pi$, 
% we shall use the usual euclidian distance $|z-z'|^2 = |x-x'|^2 + |e-e'|^2$. 
% This even if the state variable is discrete. In fact it will allow to use standard existing result, expressed only for the usual euclidian distance.

% We also define a bounded version of that euclidian distance
% \begin{equation} \label{eq:bounded_dist}
% |z-z'|_b := \min(1, |z-z'|).
% \end{equation}
% The interest of that distance is that it saturates when the two states are distinct :  $|z-z'|_b =1$  its maximal value when $e \neq e'$. So that when estimating the evolution of that quantity between two random process along time, it is not necessary to take care of the case where the state are distinct (at least when you only need a bound by above). 
% }

\medskip
\paragraph{\bf Wasserstein distance} We now introduce the Wasserstein distance together with its variants which are the main tools in the estimates for the convergence of the empirical measures throughout this paper.

For $p \geq 1$, denote by $\mathcal{P}_p(\Pi)$ the set of probability measures $\mu$ on $\Pi$ with finite $p$-th order moment, i.e. satisfying $\int_{\Pi}|z|^p\mu(dx,de)<\infty$. For $\mu$ and $\nu$ in $\mathcal{P}_p(\Pi)$, we define the $p$-Wasserstein distance by
\[
W_p^p(\mu,\nu)=\inf\left\{ \int_{\Pi\times \Pi} |z-z'|^p m (dz,dz'):\; m\in \mathcal{P}_p(\Pi\times \Pi)\text{ with marginals } \mu \text{ and }\nu \right\}.
\]
It can also be reformulated using random variables as follows,
\[
W_p^p(\mu,\nu)=\inf\left\{ \E |Z-Z'|^p  \vert \; (Z,Z') \; r.v. \text{ s.t. } \mathcal L(Z) = \mu \text{ and }\mathcal L(Z') = \nu \right\}.
\]

This distance induces the topology of weak convergence of measures together with the convergence of all moments of order up to $p$. When $p=1$, the 1-Wasserstein distance is equivalent to the bounded Lipschitz distance
\[
W_1(\mu,\nu)=\sup\left\{\int_{\Pi}\phi(x)\mu(dx)-\int_{\Pi}\phi(x)\nu(dx):\; \phi : \Pi \to \R \text{ with } Lip(\phi)\leq 1\right\}.
\]

Traditionally, to deal with the diffusion terms in the equation, it seems convenient to use the $2$-Wasserstein distance. But in this work, we have some compelling reasons to use the $1$-Wasserstein distance rather than other order $p$ (with $p>1$) due to the estimates for the jumps on the epidemiological states.

%\begin{proposition}
%		For $p\geq 1$, given $(\Omega, \F, (\F_t)_{t\in [0,T]}, \Pb)$ and $X_t$ is a martingale with respect to $\F_t$. There exists some constant $C_p>0$ such that
%		$$\E\left[\sup_{0 \leq s \leq t}|X_s|^p|\F^0\right]\leq C_p\E\left[[X]_t^{p/2}|\F^0\right]$$
%	\end{proposition}
%	The above is a conditional version of Burkholder-Davis-Gundy inequality. We will use it many times in the next proofs.

\medskip
\subsection{Conditionnal propagation of chaos}
In our setting, individuals are always linked through the common noise so there is no way to expect any full independence in the limit. However, we expect the conditional independence, given common noise. That is why we need a conditional version of chaoticity. We will reconsider some definitions and basic properties in the classical theory of mean field analysis of interacting pacticle systems, then adapt them to a frame with the presence of a common environmental noise. Before stating properly the conditional chaoticity, we need a definition of conditional exchangeability. First, recall that a sequence of random variables is said to be exchangeable if its joint distribution is invariant under any finite permutation.

\begin{definition}\label{def:cond-exchangeability}
    A sequence of $N$ random variables $\left((Z^{i,N})_{i\leq N}\right)_{N\in \mathbb{N}}$ is said to be exchangeable conditionally on $\mathcal F^0$ if, for any permutation $\sigma$ of $\{1,\dots,N\}$,
    \[\Pb^0-\text{a.s. }, \quad
    \mathcal L \Bigl( \bigl(Z^{\sigma(i),N}\bigr)_{i \le N} \vert \mathcal F^0 \Bigr)
    = \mathcal L \Bigl( \bigl(Z^{i,N}\bigr)_{i \le N} \vert \mathcal F^0\Bigr).
    \]
\end{definition}

\begin{definition}\label{def:cond-chaos}
	Let $\mu$ be a $\mathcal F^0$-measurable random probability measure on $\Pi$ or on $C([0,T], \Pi)$. The sequence $\left((Z^{i,N})_{i\leq N}\right)_{N\in \mathbb{N}}$ of conditionally exchangeable random variables is called $\mu$-chaotic if it satisfies any of three following equivalent properties :
\begin{itemize}
\item[$(i)$]
$ \displaystyle 
\Pb^0-\text{a.s. }, \; \forall \,k \in \mathbb{N},   \quad \mathcal{L}\bigl((Z^{1,N},\dots,Z^{k,N}) \vert \mathcal F^0 \bigr)\xrightarrow[N \to \infty]{\mathcal L}  \mu^{\otimes k}$,
\item[$(ii)$]
$ \displaystyle \Pb^0-\text{a.s. },  \quad \mathcal{L}\bigl((Z^{1,N},Z^{2,N}) \vert \mathcal F^0 \bigr) \xrightarrow[N \to \infty]{\mathcal L} \mu^{\otimes 2}$,
\item[$(iii)$] 
$ \displaystyle \Pb^0-\text{a.s. },  \quad 
\mu^N = \frac1N \sum_{i=1}^N \delta_{Z^{i,N}} \xrightarrow[N \to \infty]{\mathcal L} \mu.$
\end{itemize}
\end{definition}

The above definition is an adaptation of the classical one from~\cite{Sznitman} (See also~\cite{Hauray-Mischler} for a quantified version). Apparently, the conditional exchangeability is the key property in that definition.  We will rely on this assumption in order to prove the equivalence of the three properties in Definition \ref{def:cond-chaos}. The proof is in the Appendix~\ref{app:cond_chaos}.

\medskip
Under the assumption that initially the individuals are conditional exchangeable, and thanks to the unique solvability of system \eqref{system}, the conditional exchangeability property persists in time. Indeed, we state that property via the following lemma.

\begin{lemma}\label{lemma:cond-exchangeability}
The conditional exchangeability assumption holds trajectorially, i.e.
\[\Pb^0-\text{a.s.}, \; 
\forall \, \sigma  \in \Sigma^N, \quad 
\mathcal L \Bigl( \bigl(Z^{\sigma(i),N}\bigr)_{i \le N} \Bigr)
= \mathcal L \Bigl( \bigl(Z^{i,N}\bigr)_{i \le N} \Bigr).
\]
where $\Sigma^N$ denote the set of permutations of $\{1,\dots,N\}$.

In particular, at any time $t$, for all $\phi \in C_{b}\left(\Pi^{N}\right)$,
\[\mathbb{E}\left[\phi\left(Z_{t}^{1, N}, \ldots, Z_{t}^{N, N}\right) \vert \mathcal{F}_{t}^{0}\right]=\mathbb{E}\left[\phi\left(Z_{t}^{\sigma(1), N}, \ldots, Z_{t}^{\sigma(N), N}\right) \vert \mathcal{F}_{t}^{0}\right].\]
\end{lemma}

The proof of Lemma \ref{lemma:cond-exchangeability} can be found in the Appendix~\ref{app:cond_exchangeability}.

\subsection{Main results}
Before introducing the exact limit of the stochastic interacting system \eqref{system}, let us present a heuristic derivation of how the limit should look like. 

The individuals are driven by independent diffusions, a drift term, jump terms with mean-field rate, and a common diffusion. So they will soon be strongly correlated by the action of the common diffusion. A possible limit towards a mean-field equation is necessarily random, due to that common noise.

We can reasonnabily assume that the individuals will remain almost conditionally independent along the dynamics. In fact, the influence of one individual on another one in the mean field terms $\frac1N \sum_j V(Y^{i,N}_t,Y^{j,N}_t)$ and $\frac1N \sum_j K(X^{i,N}_t,X^{j,N}_t)$ is of order $\frac1N$ and we can expect that it will remain small. So we can reasonably expect that in the limit $N \to \infty$, and conditionnally to the common noise, the individual $i$ becomes independent from the other ones, and follows an equation similar to~\eqref{system}, with the mean field terms replaced by conditionnal (w.r.t. the common noise) expectations.

In the sequel, we will show that this formal argument leads at least to the right limit. The main result of this paper show that the sequence $(X^{1,N},\ldots,X^{N,N})_{N \in \N}$ is in fact $X$-chaotic, where $X$ is a process solution to the following non-linear jumping SDE (or McKean-Vlasov jump diffusion process),
\begin{equation}\label{limit_SDE}
\begin{cases}
    dX_t &= V_{\mu_t}(X_t,E_t) dt+\sigma(X_t,E_t)dB_t+\sigma_0(X_t,E_t)dB^0_t,\\
    E_t &= E_0 + P \left( \displaystyle\int_0^t   \Bigl(K_{\mu_s}(X_{s}) \1_{E_{s^-}=0} + \gamma \1_{E_{s^-}=1} \Bigr) \,ds \right),\\
    \mu_t &= \mathcal{L}\left((X_t, E_t)\vert \F^0_t\right)
	\end{cases}
\end{equation}
where $\mu_t$ is the conditional law of $(X_t,E_t)$ given $\F^0_t$,
%In both equations in system \eqref{system}, there are terms involving of empirical averaging over the other individuals. 
%As long as the particles are almost independent, is reasonnable to expect that when the number of individuals is large, these terms converges towards their expected value. 
%Now let us discuss the form of the limit of the empirical measures. 
%Supposing that $\mu_t$ the limit of the empirical measure $\mu^N_t$ exists, it should be of the form
%Here the conditioning reflects the dependencies between the individuals and under the information generated by the Brownian motion $B^0_t$, they become conditionally independent.
%\medskip
% \paragraph{\bf The mean field terms}
and the notations $V_{\mu} : \Pi\to \Rd$ and $K_{\mu} : \R^d \to \R^d$ are defined as follows,
\begin{equation}\label{def:mean-field-terms}
\begin{aligned} 
	V_{\mu}(x,e)&=\int_{\Pi} V\bigl(x,e,y,f\bigr)\mu(dy,df),\\
	K_{\mu}(x)&=\int_{\Pi} \1_1(e)K(x,y)\mu(dy,de).
\end{aligned}
\end{equation}

%In the following, we will show that the stochastic dynamic system \eqref{system} converge to a nonlinear jumping SDE on $(\Omega,\F,(\F_t)_{t\geq 0},\Pb)$,

\medskip

First, we state a theorem about the well-posedness of this stochastic system. It is in fact necessary to have a well-posed limit SDE if we want to state a result of convergence when the number of individuals tends to infinity.

\begin{theorem}\label{exist_unique_thm}
(i) {\it(Strong existence)} Given the filtered probability space $(\Omega,\F,(\F_t)_{t\geq 0},\Pb)$. Under the Lipschitz continuous assumptions on the  kernels $K,V$ and the coefficients $\sigma,\sigma_0$, there exists a solution $Z_t:=(X_t,E_t)$ to the SDE \eqref{limit_SDE} associated to any given $\F_0$-measurable initial condition $Z_0=(X_0,E_0)$.

(ii) {\it(Strong stability)} If $(X^1_t,E^1_t),\,(X^2_t,E^2_t)$ are two solutions to \eqref{limit_SDE} built on the same probability space with the same driving noises, then for any $t\geq 0$, one has the following stability estimate
\begin{equation*}
\E \left[ \sup_{0 \leq s \leq t} \left|(X^1_s,E^1_s) - (X^2_s,E^2_s) \right| \right]\leq C(t)\E \Big[ \left|(X^1_0,E^1_0) - (X^2_0,E^2_0) \right| \Big].
\end{equation*}
\end{theorem}

The derivation of macroscopic equations from microscopic models is a classical topic in the study of interacting pacticle systems in statistical physics. For this epidemic interacting system, we introduce the Kolmogorov forward equation associated to the nonlinear SDE~\eqref{limit_SDE}, which is the following nonlinear SPDE, where for each $t\ge 0$, $\mu_t$ is a random probability on $\Pi$:
\begin{equation}\label{spde}
\begin{aligned}
d\mut    &= -\dx\cdot  \left(V_{\mu_t}\mut\right)dt + 
		 \frac{1}{2}\tr\left[\nabla^2_{xx} \big((\sigma\sigma^T+\sigma_0\sigma^T_0 ) \mut\big)\right]dt
		-\dx\cdot\left(\sigma_0\mu_t \right)d\Bot
		+\big(\Gamma^I_t+\Gamma^R_t\big)dt,
	\end{aligned}
\end{equation}
where
\begin{equation*}
\begin{aligned}
          \Gamma^I_t(x,e) & = K_{\mu_t}(x) \mut(x,0) \bigl(\1_{e=1} - \1_{e=0} \bigr), \\
    	 \Gamma^R_t(x,e) & =\gamma \mut(x,1) \bigl(\1_{e=2} - \1_{e=1} \bigr).
	\end{aligned}
\end{equation*}

In fact, looking in details at the above equation, we recognize the first term in the r.h.s. as the drift term created by the drift coefficient in~\eqref{limit_SDE}, the second as the diffusion term linked to the individual Brownian motions, the third (with a part of the second one) is the random drift term linked to the common Brownian motion, and the last two $\Gamma$ terms are related to the jumps between the three different states. We will justify this relation rigorously by the following proposition.
\begin{proposition} \label{prop:NLSDE_to_SPDE}
Assume that $X$ is a solution of equation~\eqref{limit_SDE}, again under the assumption of Theorem~\ref{exist_unique_thm}. Then the collection $(\mu_t)_{t\geq 0}$ of its marginal laws solves the above SPDE~\eqref{spde}.
\end{proposition}

In the applications, it is worth to show explicitly the evolution of the empirical measure on each state of the epidemic. We can represent the function $V:\big(\R^d\times\{0,1,2\})^2\to\R^d$ as a matrix $\big(V_{i,j}\big)_{i,j\in\{0,1,2\}}$ where its entries are functions defined on $\R^d\times\R^d$, 
and introduce the matrix $A[\mu]$ which describes for the epidemiological interactions,
\[
A[\mu] :=
\begin{pmatrix}
- K_{\mu}  & 0 & 0 \\
K_{\mu} & - \gamma & 0 \\
0 & \gamma & 0
\end{pmatrix},
\]
If we write $\mu_t$ as
\[
\mut = \,^t\bigl( \mut \1_{e=0}, \mut \1_{e=1}, \mut \1_{e=2}\bigr) = 
\,^t\bigl( \mu^S_t, \mu^I_t, \mu^R_t \bigr),
\]
then the above equation~\eqref{spde} rewrites
\begin{equation}\label{spde-matrix}
	d\mut= -\dx\cdot\left(V_{\mu_t}\mut\right)dt + 
	\frac{1}{2}\tr\left[\nabla^2_{xx} \big((\sigma\sigma^T+\sigma_0\sigma^T_0 ) \mut\big)\right]dt
	-\dx\cdot\left(\mut\sigma_0 \right)d\Bot
	+ A[\mu_t] \mut dt,
\end{equation}
where the differential operators act on each component of $\mut$.

\medskip
\paragraph{\bf The results about the conditional propagation of chaos}

Back to the epidemic dynamic system \eqref{system}, our goal is to prove that this system converges in an appropriate sense to the McKean-Vlasov SDE \eqref{limit_SDE}. For that pupose, we will need to extend the classical results of propagation of chaos to the conditional one.

Following a standard strategy, we will couple the interacting system to
an auxiliary system, made of $N$ i.i.d. copies of the limit SDE, constructed on the same probability space with the same initial conditions (or at least well coupled initial conditions) and the same driving random process as for the interacting system.

Then, the estimates used in the proof of the stability result of Theorem~\ref{exist_unique_thm} together with some previously known bound (in Wasserstein sense) on the speed of convergence in the empirical law of large number~\cite{Fournier-Guillin} adapted to the case of conditionally independent random variables in Appendix~\ref{app:FG_new} will allow us to establish the following quantitative result about the conditionnal propagation of chaos for the interacting system of $N$ individuals.

On the same probability space $(\Omega,\F,(\F_t)_{t\geq 0},\Pb)$ introduced before, for each $i\in \{1,\dots, N\}$, we denote by $\bar{Z}^i=(\bar{X}^i_t,\bar{E}^i_t)$ the strong solution to the McKean-Vlasov jump-diffusion SDE similar to \eqref{limit_SDE} but for (independent) driving random processes $B^i,P^i$ instead of $B,P$ and subjects to the initial conditions $\left\{(\bar{X}^i_0,\bar{E}^i_0)\right\}_{1\leq i \leq N}$. Note that by theorem \eqref{exist_unique_thm}, for each $i$, that SDE has a unique solution and obviously the conditional law $\mathcal{L}\left((\bar{X}^i_t, \bar{E}^i_t)\big|\F^0_t\right)$ is the probability measure $\mut$ which solves \eqref{spde}.

\begin{theorem}\label{Prop_control_couplings}
	Assume that the initial conditions $\{(X^{i,N}_0,E^{i,N}_0)\},\{(\bar{X}^i_0,\bar{E}^i_0)\}, i=1,\dots,N$ are conditionally exchangeable, and have a finite moment of order $q$ for some $q>2$. Then, there exists a constant $C(t)$ such that, for all $i \in \{1,\dots,N\}$,
	\begin{equation*}
		\E \left[ \sup_{0 \leq s \leq t} \left|X^{i,N}_s-\bar{X}^i_s\right|+|E^{i,N}_s-\bar{E}^i_s|\right]\leq C(t)\left(\E \left[|X^{i,N}_0-\bar{X}^i_0|+|E^{i,N}_0-\bar{E}^i_0|\right]+\alpha_d(N)\right)
	\end{equation*}
	where $\alpha_d(N)=N^{-1/2}$ if $d=1$, $\alpha_d(N)=N^{-1/2}\log(N)$ if $d=2$, and $\alpha_d(N)=N^{-1/d}$ if $d\geq 3$.
\end{theorem}

This result is trajectorial. If we only care about the evolution of the empirical measure, then we could get a simpler result expressed in terms of empirical measures only.

\begin{corollary}\label{rate_of_convergence}
	Suppose that the random measure $\mu_0$ has a finite moment of order $q$ for some $q>2$, and that $K$,$V$,$\sigma$,$\sigma_0$ are Lipschitz continuous (as in Theorem~\ref{exist_unique_thm}). Then for all $t >0$, there exists a constant $C(t)$ such that
	\begin{align*}
		%\sup_{0 \leq s \leq t}
		\E\left[W_1\left(\mun,\mu_t\right)\right]\leq C(t)\left(\E\left[W_1\left(\mu^N_0,\mu_0\right)\right]+\alpha_d(N)\right)
	\end{align*}
	where $\alpha_d(N)=N^{-1/2}$ if $d=1$, $\alpha_d(N)=N^{-1/2}\log(N)$ if $d=2$, and $\alpha_d(N)=N^{-1/d}$ if $d\geq 3$.
\end{corollary}

%%%%%%%%%%%%%%%%%%%%%%%%%

\section{Proof of theorem \ref{exist_unique_thm} : Well-posedness of the nonlinear SDE}
\label{sec:proof_NLSDE}

\begin{proof}[Proof of ii)]
    We first establish a stability estimate for the solutions of the system \eqref{limit_SDE} which is the main estimate for the rest of this paper.
	
	Let $(X^1_t,E^1_t),\,(X^2_t,E^2_t)$ be two solutions to the SDE \eqref{limit_SDE} and $\mu^1_t, \,\mu^2_t$ respectively be conditional law of $(X^1_t,E^1_t),\,(X^2_t,E^2_t)$ given $\F^0_t$.
	
	We have
	\begin{align*}
	X^1_t-X^2_t=& X^1_0-X^2_0+\displaystyle\int_0^t\left( V_{\mu^1_s}(X^1_s,E^1_s)-V_{\mu^2_s}(X^2_s,E^2_s)\right)ds\\&+\displaystyle\int_0^t\left(\sigma(X^1_s,E^1_s)-\sigma(X^2_s,E^2_s)\right)dB_s+\displaystyle\int_0^t\left(\sigma_0(X^1_s,E^1_s)-\sigma_0(X^2_s,E^2_s)\right)dB^0_s,\\
	E^1_t-E^2_t =& E^1_0-E^2_0+\displaystyle\int_{[0,t]\times\R_+}\1_{\big\{u\leq K_{\mu^1_s}(X^1_s)\1_0(E^1_s)+\gamma \1_1(E^1_s)\big\}}\Q(ds,du)\\&-\displaystyle\int_{[0,t]\times\R_+}\1_{\big\{u\leq K_{\mu^2_s}(X^2_s)\1_0(E^2_s)+\gamma \1_1(E^2_s)\big\}}\Q(ds,du),
	\end{align*}
	which leads to
	\begin{align*}
	\sup_{0\leq s\leq t}\big|X^1_s-X^2_s\big|
	\leq|X^1_0-X^2_0|&+ \underset{0\leq s\leq t}{\sup}\left|\int_{0}^{s}\left(V_{\mu^1_r}(X^1_r,E^1_r)-V_{\mu^2_r}(X^2_r,E^2_r)\right)dr\right|\\
	&+\underset{0\leq s\leq t}{\sup}\left | \int_{0}^{s}(\sigma(X^1_r,E^1_r)-\sigma(X^2_r,E^2_r))dB_r\right |\\
	&+\underset{0\leq s\leq t}{\sup}\left | \int_{0}^{s}(\sigma_0(X^1_r,E^1_r)-\sigma_0(X^2_r,E^2_r))dB^0_r\right |
	\end{align*}
	and
	\begin{align*}
	\underset{0\leq s\leq t}{\sup}|E^1_s-E^2_s|\leq& \left|E^1_0-E^2_0\right|+ \underset{0\leq s\leq t}{\sup}\int_{[0,s]\times\R_+}\left|\1_{\left\{u\leq \alpha^1_r\right\}}-\1_{\left\{u\leq \alpha^2_r\right\}}\right|\Q(dr,du)
	\end{align*}
	where
	\begin{align*}
	    \alpha^1_r:=& \left(K_{\mu^1_r}(X^1_r)\1_0(E^1_r)+\gamma \1_1(E^1_r)\right),\\
	    \alpha^2_r:=& \left(K_{\mu^2_r}(X^2_r)\1_0(E^2_r)+\gamma \1_1(E^2_r)\right).
	\end{align*}
	
	We now take the expectation on both sides and our goal is to establish locally uniform in time bounds for the displacement and the jump between the pairs $(X^1_t,E^1_t),(X^2_t,E^2_t)$. For any $s\leq t$, we denote by $m_s$ an arbitrary coupling of $\left(\mu^1_s,\mu^2_s\right)$ on $\Pi\times\Pi$.
	
	\medskip
	\paragraph{\bf Step 1.}
	First, we will give a bound for the drift term in the displacement.
	\begin{align*}
	\E\Bigg[\sup_{0\leq s\leq t}&\left|\int_{0}^{s}\left(V_{\mu^1_r}(X^1_r,E^1_r)-V_{\mu^2_r}(X^2_r,E^2_r)\right)dr\right|\Bigg]\\
	\leq &\E\left[\int_{0}^{t}\left|V_{\mu^1_s}(X^1_s,E^1_s)-V_{\mu^2_s}(X^2_s,E^2_s)\right|ds\right]\\
	=&\E\left[\int_0^t\left|\int_{\Pi\times\Pi}\Bigl(V\bigl(X^1_s,E^1_s,x,e\bigr)-V\bigl(X^2_s,E^2_s,y,f\bigr)\Bigr)m_s(dx,de,dy,df)\right|ds\right]\\
	\leq &L_V\int_{0}^{t}\E\left[|X^1_s-X^2_s|+|E^1_s-E^2_s|+\int_{\Pi\times\Pi}\bigl(|x-y|+|e-f|\bigr)m_s(dx,de,dy,df)\right]ds.
	\end{align*}
% 	Since the choice of the coupling $m_s$ is arbitrary, we obtain the following bound due to the definition of the $1$-Wassertein metric
% 	\begin{align*}
% 	    \E\Bigg[\sup_{0\leq s\leq t}&\left|\int_{0}^{s}\left(V_{\mu^1_r}(X^1_r,E^1_r)-V_{\mu^2_r}(X^2_r,E^2_r)\right)dr\right|\Bigg]\\
% 	    &\leq L_V\int_{0}^{t}\E\left[|X^1_s-X^2_s|+|E^1_s-E^2_s|+W_1(\mu_s^1,\mu_s^2)\right].
% 	\end{align*}
	To treat the diffusion terms, we apply the Burkholder–Davis–Gundy inequality and use the Lipschitz continuity of the coefficients $\sigma, \sigma_0$
	\begin{align*}
	\E&\Bigg[\underset{0\leq s\leq t}{\sup}\left | \int_{0}^{s}\bigl(\sigma(X^1_r,E^1_r)-\sigma(X^2_r,E^2_r)\bigr)dB_r\right |+\underset{0\leq s\leq t}{\sup}\left | \int_{0}^{s}\bigl(\sigma_0(X^1_r,E^1_r)-\sigma_0(X^2_r,E^2_r)\bigr)dB^0_r\right |\Bigg]\\
	&\leq C_1(L_{\sigma}+L_{\sigma_0})\E\left[\left(\int_{0}^{t}\bigl(|X^1_s-X^2_s|+|E^1_s-E^2_s|\bigr)^2ds\right)^{1/2}\right].
	\end{align*}
	\begin{remark}\label{rmk_a_direct_estimate}
	At this point, a typical difficulty arises where there is no way to remove the powers or simplify the terms. One possibility is to take the supremum in time inside the integral and bound it first on a small enough time interval as in \cite{Coghi-Flandoli} or \cite{Erny-Locherbach-Loukianova}. Once we have an estimate on a small time interval, the proof is done by iterating over the remain intervals. Nevertheless, in this paper we will make a slight modification by controlling the quantities inside the integral in order to obtain directly an estimate on an arbitrary time interval.
	\end{remark}

	Indeed, one has
	\begin{align*}
	    C_1&(L_{\sigma}+L_{\sigma_0})\E\left[\left(\int_{0}^{t}\bigl(|X^1_s-X^2_s|+|E^1_s-E^2_s|\bigr)^2ds\right)^{1/2}\right]\\
	    \leq& \E\left[C_1(L_{\sigma}+L_{\sigma_0})\left(\sup_{0\leq s\leq t}\bigl(|X^1_s-X^2_s|+|E^1_s-E^2_s|\bigr)\right)^{1/2}\left(\int_{0}^{t}\bigl(|X^1_s-X^2_s|+|E^1_s-E^2_s|\bigr)ds\right)^{1/2}\right]\\
	    \leq& \frac{1}{2}\E\left[\sup_{0\leq s\leq t}\bigl(|X^1_s-X^2_s|+|E^1_s-E^2_s|\bigr)\right]+\frac{C_1^2(L_{\sigma}+L_{\sigma_0})^2}{2}\int_{0}^{t}\E\left[|X^1_s-X^2_s|+|E^1_s-E^2_s|\right]ds.
	\end{align*}
	Now summing up all the above terms on both sides, we obtain an appropriate estimate for the displacements
	\begin{equation}\label{displacement_estimate}
	\begin{aligned}
	\E\left[\sup_{0\leq s\leq t}\left|X^1_s-X^2_s\right|\right]\leq &\E\left[|X^1_0-X^2_0|\right]+\frac{1}{2}\E\left[\sup_{0\leq s\leq t}\bigl(|X^1_s-X^2_s|+|E^1_s-E^2_s|\bigr)\right]\\
	&+\left(L_V+\frac{C_1^2(L_{\sigma}+L_{\sigma_0})^2}{2}\right)\int_{0}^{t}\E\left[|X^1_s-X^2_s|+|E^1_s-E^2_s|\right]ds\\
	&+L_V\int_{0}^{t}\E\left[\int_{\Pi\times\Pi}\bigl(|x-y|+|e-f|\bigr)m_s(dx,de,dy,df)\right]ds.
	\end{aligned}
	\end{equation}

\paragraph{\bf Step 2.}
For the Poisson part with the jumps on the epidemiological states, it could be easier to deal with the supremum in time. Indeed, notice that
	\[\bigl|\1_{u \le \alpha^1_r} - \1_{u \le \alpha^2_r} \bigr| = \1_{u \in [\alpha^1_r \wedge \alpha^2_r, \alpha^1_r \vee \alpha^2_r]}\]
	and
	\[\alpha^1_r \vee \alpha^2_r - \alpha^1_r \wedge \alpha_r^2 = |\alpha^1_r - \alpha^2_r|\]
	Using the above identities, we obtain the following estimate
	\begin{align*}
	\E\bigg[\sup_{0\leq s\leq t}&\int_{[0,s]\times\R_+}\left|\1_{\left\{u\leq \alpha^1_r\right\}}-\1_{\left\{u\leq \alpha^2_r\right\}}\right|\Q(dr,du)\bigg]\\
	\leq& \E\left[\int_{[0,t]\times\R_+}\1_{\left\{u\in\left[\alpha^1_s\wedge\alpha^2_s,\alpha^1_s\vee\alpha^2_s\right]\right\}}\Q(ds,du)\right]\\
	=&\E\left[\int_0^t\Bigl|\left(K_{\mu^2_r}(X^2_r)\1_0(E^2_r)+\gamma \1_1(E^2_r)\right)-\left(K_{\mu^1_r}(X^1_r)\1_0(E^1_r)+\gamma \1_1(E^1_r)\right)\Bigr|dr\right]\\
	\leq &\E\left[\int_0^t\left| K_{\mu^1_s}(X^1_s)\1_0(E^1_s)-K_{\mu^2_s}(X^2_s)\1_0(E^2_s)\right|ds\right]+\gamma\E\left[\int_0^t\left|\1_1(E^1_s)-\1_1(E^2_s)\right|ds\right]\\
	=&I_1+I_2.
	\end{align*}
	\begin{remark}\label{rmk_using_W1}
	In the first inequality, we can see that no matter the initial degree of the integrand is, we always end with $\1_{u \in [\alpha^1_r \wedge \alpha^2_r, \alpha^1_r \vee \alpha^2_r]}$. This is usual when dealing with the Poisson processes. Hence, the following estimates will keep that of order $1$ and it is impossible to recover the initial order to get an inequality in type of the Gronwall lemma in the end. Therefore, this crucial point forces us to use the $1$-Wasserstein distance instead of the other orders in the proof, even though it may less simple to treat the diffusion terms as above.
	\end{remark}

    Since the intensity of the jump processes also depends on the mean field terms, we obtain  To treat the first term, we will use the same coupling $m_s$ as introduced previously. One has,
	\begin{align*}
	I_1=&\E\left[\int_0^t\left| K_{\mu^1_s}(X^1_s)\1_0(E^1_s)-K_{\mu^2_s}(X^2_s)\1_0(E^2_s)\right|ds\right]\\ =&\E\left[\int_0^t\left|\int_{\Pi}K(X^1_s,x)\1_0(E^1_s)\1_1(e)\mu^1_s(dx,de)-\int_{\Pi}K(X^2_s,y)\1_0(E^2_s)\1_1(f)\mu^2_s(dy,df)\right|ds\right]\\
	\leq&\E\left[\int_0^t\left|\int_{\Pi\times\Pi}\Bigl(K(X^1_s,x)-K(X^2_s,y)\Bigr)\1_0(E^1_s)\1_1(e)m_s(dx,de,dy,df)\right|ds\right]\\
	&+\E\left[\int_0^t\left|\int_{\Pi\times\Pi}K(X^2_s,y)\Bigl(\1_0(E^1_s)\1_1(e)-\1_0(E^2_s)\1_1(f)\Bigr)m_s(dx,de,dy,df)\right|ds\right]\\
	\leq& L_K\int_{0}^{t}\E\left[|X^1_s-X^2_s|+\int_{\Pi\times\Pi}|x-y|m_s(dx,de,dy,df)\right]ds\\
	&+\|K\|_{\infty}\int_0^t\E\left[|E^1_s-E^2_s|+\int_{\Pi\times\Pi}|e-f|m_s(dx,de,dy,df)\right]ds,
	\end{align*}
	where to obtain the last inequality, we used the fact that $|\1_k(e)-\1_k(f)|\leq |e-f|$, for $k\in \{0,1,2\}$. The same observation leads to
\[I_2=\gamma\E\left[\int_0^t\left|\1_1(E^1_s)-\1_1(E^2_s)\right|ds\right]\leq \gamma\int_0^t\E\left[\left|E^1_s-E^2_s\right|\right]ds.
\]
	Summing up the above terms, we obtain the following estimate for the changes between the epidemiological states
	\begin{align}\label{jump_estimate}
	\begin{split}
	\E\left[\underset{0\leq s\leq t}{\sup}|E^1_s-E^2_s|\right]\leq &\E\left|E^1_0-E^2_0\right|+L_K\int_{0}^{t}\E\left[|X^1_s-X^2_s|\right]ds\\
	&+\bigl(\gamma+\|K\|_{\infty}\bigr)\int_{0}^{t}\E\left[|E^1_s-E^2_s|\right]ds\\
	&+ L_K\int_{0}^{t}\E\left[\int_{\Pi\times\Pi}|x-y|m_s(dx,de,dy,df)\right]ds\\
	&+ \|K\|_{\infty}\int_0^t\E\left[\int_{\Pi\times\Pi}|e-f|m_s(dx,de,dy,df)\right]ds.
	\end{split}	
	\end{align}
	
	\paragraph{\bf Step 3.}
	We combine \eqref{displacement_estimate} and \eqref{jump_estimate} to finally deduce
	\begin{align*}
	\frac12\E \bigg[ \sup_{0 \leq s \leq t} & \left|X^1_s-X^2_s\right|+|E^1_s-E^2_s|\bigg]\\
	\leq &\E \left[|X^1_0-X^2_0|+|E^1_0-E^2_0|\right]\\
	&+\left(L_V+\frac{C_1^2(L_{\sigma}+L_{\sigma_0})^2}{2}+L_K\right)\int_{0}^{t}\E\left[|X^1_s-X^2_s|\right]ds\\
	&+\left(L_V+\frac{C_1^2(L_{\sigma}+L_{\sigma_0})^2}{2}+\gamma+\|K\|_{\infty}\right)\int_{0}^{t}\E\left[|E^1_s-E^2_s|\right]ds\\
	&+\left(L_V+L_K \right)\int_{0}^{t}\E\left[\int_{\Pi\times\Pi}|x-y|m_s(dx,de,dy,df)\right]ds\\
	&+\left(L_V+\|K\|_{\infty}\right)\int_0^t\E\left[\int_{\Pi\times\Pi}|e-f|m_s(dx,de,dy,df)\right]ds.
	\end{align*}
	Since $m_s$ is the conditional law of the pair $\big((X^1_s,E^1_s),(X^2_s,E^2_s)\big)$ given $\F^0$, we can simplify the last two terms and deduce the following bound
	\begin{align}\label{main_estimate}
		\begin{split}
			\E& \left[ \sup_{0 \leq s \leq t} \left|X^1_s-X^2_s\right|+|E^1_s-E^2_s|\right]\\
			\leq & 2\E\left[|X^1_0-X^2_0|+|E^1_0-E^2_0|\right]+C\int_{0}^{t}\E\left[|X^1_s-X^2_s|+|E^1_s-E^2_s|\right]ds.
		\end{split}
	\end{align}
	where $C=\max\Bigl\{\left(4L_V+4L_K+C_1^2(L_{\sigma}+L_{\sigma_0})^2\right), \left(4L_V+4\|K\|_{\infty}+2\gamma+C_1^2(L_{\sigma}+L_{\sigma_0})^2\right)\Bigr\}$.
	
	Finally, we use the Gronwall lemma to get the stability estimate
	\begin{equation}\label{stability}
	\E \bigg[ \sup_{0 \leq s\leq t}             \left|X^1_s-X^2_s\right|+|E^1_s-E^2_s|\bigg]\\
    \leq 2e^{Ct}\E \left[|X^1_0-X^2_0|+|E^1_0-E^2_0|\right].
    \end{equation}

    \medskip
    
    \paragraph{\textit{Proof of i).}}
	 To prove the existence of a solution to equation~\eqref{limit_SDE}, we will use a fixed point theorem. We consider $\Theta$ the space consisting of  $\Pi$-valued $\F$-progressively measurable \textit{càdlàg} processes $Z=(Z_t)_{t \in [0,T]}$ on $[0,T]$ equipped with the norm
	$$\|Z\|_{\Theta} = \E\Big[\sup_{0 \leq t \leq T}|Z_t|\Big]< \infty.$$
	Notice that $Z=(X,E)$ is a mixed process, where $X$ is continuous and $E$ is a \textit{càdlàg} process.
    With the norm defined above, $\left(\Theta,\|\cdot\|_{\Theta}\right)$ is a Banach space.
	
	Now we fix $Z^1=(X^1,E^1)$ and $Z^2=(X^2,E^2)$ in $\Theta$ and define a mapping $\Phi:\Theta\mapsto\Theta$ which maps $(X,E)\mapsto (Y,F)$ respectively via the following stochastic integrations:
	\begin{align*}
	\forall \, t \in [0,T], \quad    Y_t=& X_0+\int_0^t V_{\mu_s}(X_s,E_s) ds+\int_0^t\sigma(X_s,E_s)dB_s+\int_0^t\sigma_0(X_s,E_s)dB^0_s,\\
	    F_t =& E_0 + \int_{[0,t]\times\R_+}  \1_{ \left\{ u \le \{K_{\mu_s}(X_s) \1_0(E_s) + \gamma \1_1(E_s)\right\}}\Q(ds,du).
	\end{align*}

We can show that if $(X,E) \in \Theta$, then $(Y,F)=\Phi(X,E) \in\Theta$ as well by using the Burkholder–Davis–Gundy inequality and a similar estimate as in the proof of Proposition \ref{finite_moments}.
% 	\begin{align*}
%         \E\Big[\sup_{0\leq t\leq T}|Y_t|\Big]\leq& \E\big[|X_0|\big]+\E\left[\sup_{t\leq t\leq T}\left|\int_0^t V_{\mu_s}(X_s,E_s)ds\right|\right]\\
%         &+\E\left[\sup_{t\leq t\leq T}\left|\int_0^t\sigma(X_s,E_s)dB_s\right|\right]+\E\left[\sup_{t\leq t\leq T}\left|\int_0^t\sigma_0(X_s,E_s)dB^0_s\right|\right]\\
%         \leq& \E[|X_0|]+\E\left[\int_0^T\int_{\Pi}\left| V(X_s,E_s,x,e)\right|\mu_s(dx,de)ds\right]\\
%         &+ C_{BDG}\E\left[\left(\int_0^T\left|\sigma(X_s,E_s)\right|^2ds\right)^{1/2}\right]+C_{BDG}\E\left[\left(\int_0^T\left|\sigma_0(X_s,E_s)\right|^2ds\right)^{1/2}\right]\\
%         \leq& \E[|X_0|]+C(V,\sigma,\sigma_0)T
%     \end{align*}

If for $i=1,2$, $(Y^i,F^i)=\Phi(X^i,E^i)$, then by the same calculations those lead to the main estimate \eqref{main_estimate}, we obtain
\begin{equation}\label{contraction_bound}
    \begin{aligned}
	\E &\left[ \sup_{0 \leq s \leq t} |Y^1_s-Y^2_s|+|F^1_s-F^2_s|\right]\\
	&\leq \frac{1}{2} \E\left[\sup_{0\le s\le t}|X^1_s-X^2_s|+|E^1_s-E^2_s|\right]+C\int_0^t\E\big[|X^1_s-X^2_s|+|E^1_s-E^2_s|\big]ds.
    \end{aligned}
\end{equation}
In other words, one has
\begin{equation*}
    \|\Phi(Z^1) - \Phi(Z^2)\|_{\Theta}\leq \left(\frac{1}{2} +C\int_0^Tdt_1\right)\|Z^1 - Z^2\|_{\Theta}.
\end{equation*}
It follows that
\begin{align*}
	\|\Phi^2(Z^1) - \Phi^2(Z^2)\|_{\Theta}\leq &\frac{1}{2}\|\Phi(Z^1) - \Phi(Z^2)\|_{\Theta} +C\int_0^T\E\big[|\Phi(Z^1)_{t_1} - \Phi(Z^2)_{t_1}|\big]dt_1\\
	\leq &\frac{1}{2}\left(\frac12+C\int_0^Tdt_1\right)\|Z^1 - Z^2\|_{\Theta}+C\int_0^T \left(\frac12+C\int_0^{t_1}dt_2\right)\|Z^1 - Z^2\|_{\Theta}dt_1\\
	= & \left[\left(\frac{1}{2}\right)^2+ 2\left(\frac12\right)C\int_0^T dt_1+C^2\int_0^T dt_1\int_0^{t_1} dt_2\right]\|Z^1 - Z^2\|_{\Theta}.
\end{align*}
and in general,
\begin{align*}
	\|\Phi^k(Z^1) - \Phi^k(Z^2)\|_{\Theta}\leq &\frac{1}{2}\|\Phi^{k-1}(Z^1) - \Phi^{k-1}(Z^2)\|_{\Theta} +C\int_0^T\E\big[|\Phi^{k-1}(Z^1)_{t_1} - \Phi^{k-1}(Z^2)_{t_1}|\big]dt_1\\
	\leq &\cdots\\
	\leq & \bigg[\left(\frac12\right)^k+\binom{k}{k-1}\left(\frac12\right)^{k-1}C\int_0^T dt_1+\dots\\
	&\hspace{2cm}+\binom{k}{1}\left(\frac12\right)C^{k-1}\int_0^Tdt_1\dots\int_0^{t_{k-2}}dt_{k-1}\\
	&\hspace{2cm}+C^k\int_0^Tdt_1\dots\int_0^{t_{k-2}}dt_{k-1}\int_0^{t_{k-1}}dt_{k}\bigg] \|Z^1 - Z^2\|_{\Theta}.
\end{align*}
We obtain the last bound by iterating k times the inequality \eqref{contraction_bound}. Therefore, we conclude that for any $k\in \N$,
\begin{align*}
	\|\Phi^k(Z^1) - \Phi^k(Z^2)\|_{\Theta}\leq & \sum_{j=0}^k\frac{k!}{j!(k-j)!}\frac{1}{2^{k-j}}\frac{(CT)^j}{j!}\|Z^1-Z^2\|_{\Theta}.
\end{align*}

In the next step, we will show that \begin{equation}\label{constant_contraction}
    \sum_{j=0}^k\frac{k!}{j!(k-j)!}\frac{1}{2^{k-j}}\frac{(CT)^j}{j!}\xrightarrow{k\to \infty} 0,
\end{equation}
and this will imply that for $k$ lagre enough, $\Phi^k$ is a contraction on $\Theta$. Therefore, $\Phi$ has a unique fixed point.

Now, to complete the proof, we will show that the above sum \eqref{constant_contraction} vanishes in the limit $k \to \infty$. We will use the following observations:

Let $M>2CT$ be large enough such that for all $j>M$, one has $j!>M^j$. First we can handle the remaining part of the above sum for large $k$,
\begin{equation*}
    \begin{aligned}
        \sum_{j=M+1}^k\frac{k!}{j!(k-j)!}\frac{1}{2^{k-j}}\frac{(CT)^j}{j!}<& \sum_{j=M+1}^k\frac{k!}{j!(k-j)!}\frac{1}{2^{k-j}}\frac{(CT)^j}{M^j}\\
        <& \sum_{j=0}^k\frac{k!}{j!(k-j)!}\frac{1}{2^{k-j}}\frac{(CT)^j}{M^j}\\
        =&\left(\frac{1}{2}+\frac{CT}{M}\right)^k\xrightarrow{k\to\infty}0.
    \end{aligned}
\end{equation*}
For the first $M+1$ terms, we bound the sum as follows
\begin{equation*}
    \begin{aligned}
        \sum_{j=0}^M\frac{k!}{j!(k-j)!}\frac{1}{2^{k-j}}\frac{(CT)^j}{j!}=& \sum_{j=0}^M\frac{k(k-1)\dots(k-j+1)}{2^k}\frac{(2CT)^j}{(j!)^2}\\
        <& \sum_{j=0}^M\frac{k^j}{2^k}\frac{(2CT)^j}{(j!)^2}\\
        <&\frac{(M+1)k^M}{2^k}\max_{0\le j\le M}\left\{\frac{(2CT)^j}{(j!)^2}\right\}\xrightarrow{k\to\infty}0.
    \end{aligned}
\end{equation*}

\end{proof}
\begin{remark}
	In the above proof, once we coupled two processes, we used the representation of the jumps by the Poisson random measures instead of the Poisson point processes. Even though the two representations are the same in law, the law of the couplings in these two cases are very different and the use of the Poisson random measures here plays an important role. With the Poisson random measures, most of the jumps of the coupling are done at the same time, which allows us to obtain the estimates locally uniformly in time, whereas the same does not hold for the Poisson point processes, which would involve the Skorohod topology instead of the uniform topology.
\end{remark}

\begin{remark}
With the bound~\eqref{contraction_bound}, we can always find an appropriate contraction. Indeed, to generalize the argument, let us introduce the following norm
    \[
    \|Z\|_{\infty,\lambda} = \sup_{0 \leq t } \left\{
    e^{-\lambda t} \E\Big[ \bigl(\sup_{s \le t} |Z_s| \bigr)\Big]  \right\}
    < \infty.
    \]

By \eqref{contraction_bound}, the mapping $\Phi$ satisfies
\[
\E \left[ \sup_{0 \leq s \leq t} |\Phi(Z)_s| \right] \le \frac12 \E \left[ \sup_{0 \leq s \leq t} |Z_s| \right] + C \int_0^t |Z_s| \,ds.
\]
Using the norm that defined above, we have
\begin{align*}
\E \left[ \sup_{0 \leq s \leq t} |\Phi(Z)_s| \right] & \le \frac12 \E \left[ \sup_{0 \leq s \leq t} |Z_s| \right] + C \int_0^t |Z_s| \,ds. \\
& \le  \| Z\|_{\infty, \lambda} \left( \frac12 e^{\lambda t} + C \int_0^t e^{\lambda s} \,ds \right) \\
&\le \| Z\|_{\infty, \lambda} \left( \frac12  + \frac C  \lambda \right) e^{\lambda t}. \\
\end{align*}
So that
\[
\| \Phi(Z)\|_{\infty, \lambda} \le 
 \left( \frac12  + \frac C  \lambda \right)
\| Z\|_{\infty, \lambda},
\]
and this is a contraction if $\lambda > 2C$.
\end{remark}
\medskip

	With the assumptions that the kernels and coefficients of the McKean-Vlasov equation are bounded, we have the propagation of moments up to a finite time $T$.
\begin{proposition}\label{finite_moments}
	Assume that the initial law $\mu_0$ has a finite moment of order $q$ with $q\geq 1$. For all $t\leq T$, the moment of order $q$ of $\mu_t$ is also finite.
\end{proposition}
\begin{proof}
	The proof is straightforward. First, notice that the state $e$ is always bounded by $2$, so we do not need to estimate any moment involving $e$. On the other hand,
    \[
    \E\left[\int_{\Pi}|x|^q\mu_t(dx,de)\right] = \E\left[\E \left[ |X_t|^q \vert \mathcal F^0 \right] \right] = \E \left[|X_t|^q\right].
    \]
    One has,
    \begin{align*}
        \E[|X_t|^q]\leq& C_q\E[|X_0|^q]+C_q\E\Bigg[\left|\int_0^t V_{\mu_s}(X_s,E_s)ds\right|^q\Bigg]\\
        &+C_q\E\Bigg[\left|\int_0^t\sigma(X_s,E_s)dB_s\right|^q\Bigg]+C_q\E\Bigg[\left|\int_0^t\sigma_0(X_s,E_s)dB^0_s\right|^q\Bigg]\\
        \leq& C_q\E[|X_0|^q]+C_q\E\left[\int_0^t\int_{\Pi}\left| V(X_s,E_s,x,e)\right|^q\mu_s(dx,de)ds\right]\\
        &+ C_q\E\left[\left(\int_0^t\left|\sigma(X_s,E_s)\right|^2ds\right)^{q/2}\right]+C_q\E\left[\left(\int_0^t\left|\sigma_0(X_s,E_s)\right|^2ds\right)^{q/2}\right]\\
        \leq& C_q\E[|X_0|^q]+C(q,V,\sigma,\sigma_0)T
    \end{align*}
    where the last inequality follows from by the boundedness of $V,\sigma,\sigma_0$.
    
\end{proof}

%%%%%%%%%%%%%%%%%%%%%%%%%%%%%%%%%%%%%%%%%%%%%%%%%%%%%
%\section{Well-posedness of the mean field limit equation}

\section{The non linear PDE}
\label{sec:nlpde}

The empirical measure of the interacting system converges to the conditional law of the unique solution of the limit SDE \eqref{limit_SDE}, which is also the solution of the stochastic PDE \eqref{spde}. We will give a proof of Proposition~\ref{prop:NLSDE_to_SPDE}, which shows the connection between the SPDE~\eqref{spde} and the nonlinear jumping SDE~\eqref{limit_SDE}. 

\begin{proof}[Proof of Proposition~\ref{prop:NLSDE_to_SPDE}]
Let $\phi\in C^2_b(\Pi)$ and use Ito's formula to expand
\begin{align*}
	\phi(X_t,&E_t)-\phi(X_0,E_0)\\ =&\int_{0}^{t}\dx\phi(X_s,E_s)dX_s+\frac{1}{2}\int_{0}^{t}\tr\left[\dxx\phi(X_s,E_s)d\left\langle X\right\rangle_s\right]\\
	&+\int_{[0,t]\times\R_+}\big(\phi(X_s,E_s)-\phi(X_s,E_{s^-})\big)\1_{\big\{u\leq K_{\mu_s}(X_s) \1_0(E_{s^-}) + \gamma \1_1(E_{s^-}) \big\}}Q(ds,du)\\
	=&I_1+I_2+I_3
\end{align*}

We now take the conditional expectation on the information generated up to time $t$ by the common noise $B^0$ for each term above and notice that the conditional law of $\phi(X_t,E_t)$ given $F^0_t$ is the random variable $\left<\mu_t,\phi\right>$.
\begin{align*}
	\E\left[I_1\vert \F^0_t\right]=&\E\left[\int_{0}^{t}\dx\phi(X_s,E_s)V_{\mu_s}(X_s,E_s)ds \;\bigg\vert \;  \F^0_t\right]+\E\left[\int_{0}^{t}\dx\phi(X_s,E_s)\sigma(X_s,E_s)d\Bis \;\bigg\vert \; \F^0_t\right]\\
	&+\E\left[\int_{0}^{t}\dx\phi(X_s,E_s)\sigma_0(X_s,E_s)dB^0_s \;\bigg\vert \;  \F^0_t\right]\\
	=&\int_{0}^{t}\int_{\Pi}\dx\phi(x,e)V_{\mu_s}(x,e)\mu_s(dx,de)ds+\int_{0}^{t}\int_{\Pi}\dx\phi(x,e)\sigma_0(x,e)\mu_s(dx,de)dB^0_s.
\end{align*}
In the last equality, the change of integrations w.r.t $\F^0_t$ follows a kind of Fubini's theorem for conditional expectation and stochastic integral. Indeed, we first consider the simple predictable processes then apply the tower property since $\F^0_t \subset \F_{\sigma(X_t)}\vee\F_{\sigma(E_t)}$, and the pull out property of conditional expectation, and at the end passing to limit the approximations by the dominated convergence theorem. That means, for $Z_t$ a $\bigl(\F^0_t\vee\F^{1,\dots,N}_t\bigr)$-predictable process, the following holds:
\begin{align*}
    \E\left[\int_0^tZ_sdB^i_s\vert \F^0_t\right]&=0,\\
    \E\left[\int_0^tZ_sdB^0_s\vert \F^0_t\right]&=\int_0^t\E\left[Z_s\vert \F^0_s\right]dB^0_s,\\
    \E\left[\int_0^tZ_sds\vert \F^0_t\right]&=\int_0^t\E\left[Z_s\vert \F^0_t\right]ds=\int_0^t\E\left[Z_s\vert \F^0_s\right]ds.
\end{align*}
where in the last inquality, we used the fact that $\F^0_t=\F^0_s\vee\sigma\left\{B^0_r-B^0_s,\;s\le r\le t\right\}$ and $\sigma\left\{B^0_r-B^0_s,\;s\le r\le t\right\}{\perp\!\!\!\perp}Z_s$.

For the quadratic variation term we also have
\begin{align*}
	\E\left[I_2\vert\F^0\right]=&\E\left[\frac{1}{2}\int_{0}^{t}\tr\left[\dxx\phi(X_s,E_s)\left(\sigma\sigma^T(X_s,E_s)+\sigma_0\sigma^T_0(X_s,E_s)\right)\right]ds \;\bigg\vert \; \F^0_t\right]\\
	=&\frac{1}{2}\int_{0}^{t}\int_{\Pi}\tr\left[\dxx\phi(x,e)\left(\sigma\sigma^T+\sigma_0\sigma^T_0\right)(x,e)\mu_s(dx,de)\right]ds.
\end{align*}
For the last term, since the counting process has a compensator
\[
\displaystyle\int_0^s\bigl(K_{\mu_r}(X_{r}) \1_{E_{r^-}=0} + \gamma \1_{E_{r^-}=1} \bigr)dr
\]
and the integrand $\left(\phi(X_s,E_s)-\phi(X_s,E_{s^-})\right)$ is left-continuous and therefore predictable, then the expectation of $I_3$ is equal to the conditional expectation of the integration w.r.t the compensator.
\begin{align*}
	\E\left[I_3\vert\F^0_t\right]=&\E\left[\int_{[0,t]\times\R_+}\big(\phi(X_s,E_s)-\phi(X_s,E_{s^-})\big)\1_{\big\{u\leq K_{\mu_s}(X_s) \1_0(E_{s^-}) + \gamma \1_1(E_{s^-}) \big\}}Q(ds,du)\vert\F^0\right]\\
	=&\E\left[\int_{0}^{t}\bigl(\phi(X_s,E_s)-\phi(X_s,E_{s^-})\bigr)\bigl(K_{\mu_s}(X_{s}) \1_0(E_{s^-}) + \gamma \1_1(E_{s^-}) \bigr)ds\vert\F^0\right]\\
	=&\int_{0}^{t}\int_{\Pi}\bigl(\phi(x,e+1)-\phi(x,e)\bigr)K_{\mu_s}(x) \1_{e=0}\mu_s(dx,de)ds\\
	&+\int_{0}^{t}\int_{\Pi}\bigl(\phi(x,e+1)-\phi(x,e)\bigr)\gamma \1_{e=1}\mu_s(dx,de)ds.
\end{align*}
Combining all the above terms, we deduce the weak form of the SPDE \eqref{spde}
\begin{equation*}
\begin{aligned}
    d\left<\phi,\mu_t\right>=&\left<V_{\mu_t}\cdot\dx\phi,\mu_t\right>dt+\frac{1}{2}\left<\tr\left[\left((\sigma\sigma^T+\sigma_0\sigma^T_0)\right)\dxx\phi\right],\mu_t\right>dt+\left<\sigma_0\cdot\dx\phi,\mu_t\right>dB^0_t\\
    &+\left<K_{\mu_t}(\1_{e=1}-\1_{e=0})\phi,\mu_t(dx,0)\right>dt+\left<\gamma (\1_{e=2}-\1_{e=1})\phi, \mu_t(dx,1)\right>.
\end{aligned}
\end{equation*}
\end{proof}

%\vvy{We will prove the uniqueness of solutions to the SPDE \eqref{spde} in the class of solutions $\mu=(\mu_t)_{t\in [0,T]}$ taking values in the space of finite non-negative measures on $\Pi$ (\mh{measurable in time, or weakly continuous...}).  

\begin{remark}

The question of uniqueness of solutions to the above nonlinear SPDE is interesting, and of course related to the uniqueness of solution to the McKean-Vlasov jump-diffusion process~\eqref{limit_SDE}. We can point out a strategy to prove the uniqueness of solutions to the nonlinear SPDE~\eqref{spde} in the class of measure-valued solutions relies on the uniqueness result proved for solution of the nonlinear SDE~\eqref{limit_SDE}, and also the uniqueness of solutions to a corresponding \emph{linear} SPDE (which is unknown so far).

We discuss here a sketch of that strategy. Indeed, we start with a  solution $\mu$ to the SPDE~\eqref{spde} and freeze the coefficients of that SPDE by setting $a(t,x)=K_{\mu_t}(x),\; b(t,x)=V_{\mu_t}(x)$. We obtain a linear SPDE, with smooth coefficients (precisely the smoothness),
\begin{equation}\label{linear_SPDE}
\begin{aligned}
d\nu_t    &= -\dx\cdot  \left(b\nu_t\right)dt + 
		 \frac{1}{2}\tr\left[\nabla^2_{xx} \big((\sigma\sigma^T+\sigma_0\sigma^T_0 ) \nu_t\big)\right]dt
		-\dx\cdot\left(\sigma_0\nu_t \right)d\Bot\\
		&+a\bigl(\1_{e=1} - \1_{e=0} \bigr)\nu_t(dx,0) dt+\gamma \bigl(\1_{e=2} - \1_{e=1} \bigr)\nu_t(dx,1) dt
	\end{aligned}
\end{equation}
Obviously $\mu$ solves the linear version \eqref{linear_SPDE} with coefficients $a, b$. If that equation has at most one solution in the class of finite measures, so that $\mu$ is the unique solution in that class to the linear SPDE~\eqref{linear_SPDE}.

%and prove the uniqueness of solutions to this linear version 
%(Could we extend the result of N.Fournier, L. Xu on the equivalence between jumping SDEs and its corresponding SPDEs to the conditional case, or only need to prove directly the uniqueness in the class of finte non-negative measures for this linear SPDE?)

We next consider an underlying the linear version of the McKean-Vlasov SDE with jumps \eqref{limit_SDE}, which conditional law (at fixed time) is a  solution to linear SPDE \eqref{linear_SPDE} :
\begin{equation}\label{linear_SDE}
\begin{cases}
    dX_t &= b(X_t,E_t) dt+\sigma(X_t,E_t)dB_t+\sigma_0(X_t,E_t)dB^0_t\\
    E_t &= E_0 + P \left( \displaystyle\int_0^t   \Bigl(a(X_{s}) \1_{E_{s^-}=0} + \gamma \1_{E_{s^-}=1} \Bigr) \,ds \right)\\
    \mu_t &= \mathcal{L}\left((X_t, E_t)\vert \F^0_t\right)
\end{cases}
\end{equation}

Let $\Omega$ be a probability space which is rich enough to define the Brownian motions $B$, $B^0$ and the Poisson point process  $P$, independent of each others, and also a r.v.  $Z_0$ which is a $\F_0$-measurable random variable such that $\mathcal{L}(Z_0)=\mu_0$.
Using standard results in the field of SDEs, we could define a process $Z:[0,T]\times\Omega\to \Pi$ that is a $\F_t$-adapted solution to linear SDE \eqref{linear_SDE} (with coefficients $a,b$). Using Itô's formula, we know that $\mathcal{L}\big(Z\vert \F^0\big)$ is also a solution to equation \eqref{linear_SPDE}. Now asuming the uniqueness of solutions to the linear SPDE \eqref{linear_SPDE} in the class of finite measures, we can conclude that for any $t \ge 0$, $\Pb^0-\text{a.s}$, $\mu_t=\mathcal{L}\big(Z_t\vert \F^0\big)$.

But this implies that the above process $Z$ is a solution to the limit jumping SDE \eqref{limit_SDE}. So, if we start from a second solution $\nu$ to the SPDE~\eqref{spde}, with the same initial condition $\nu_0=\mu_0$. Then, repeating the above argument, we can construct \emph{on the same probability space $\Omega$} a process $Z'$ solution to~\eqref{limit_SDE}, such that $\nu_t=\mathcal{L}\big(Z'_t\vert \F^0\big)$ for any $t\ge 0$.

Now by Theorem \eqref{exist_unique_thm} which states that the uniqueness of processes solution to equation \eqref{limit_SDE} constructed on the same probability space. This means that in fact $Z=Z'$ almost surely, and this implies that $\Pb^0-\text{a.s}$, $\mu =\nu$.

\end{remark}

%This means $\mu$ is a solution to the SPDE \eqref{spde} if and only if there exists an $\F_t$-adapted stochatic process $Z_t$ such that $\mu_t=\mathcal{L}\big(Z_t\vert \F^0\big)$ and the pair $(Z,\mu)$ is a solution to the limit jumping SDE \eqref{limit_SDE}.

%%%%%%%%%%%%%%%%%%%%%%%%%%%%%%%%%%%%%%
%
%
%%%%%%%%%%%%%%%%%%%%%%%%%%%%%%%%%%%%%%

\section{Law of large numbers and propagation of chaos}\label{sec:LLN}
In this section, we provide a quantitative estimates for the convergence of the empirical measure of the epidemic dynamic system to the conditional law given $\mathcal{F}^0$ of the unique solution of the mean field limit.

We now introduce the auxilary system as described before, which is made of independent copies of the limit SDE driven by independent noises $B^i$ instead of $B$, independent PRMs $Q^i$ instead of $Q$ namely
\begin{equation}\label{copies}
	\begin{cases}
		X^i_t &= X^i_0+\displaystyle\int_{0}^{t}V_{\mu_s}(X^i_s,E^i_s) ds+\int_{0}^{t}\sigma(X^i_s,E^i_s)dB^i_s+\int_{0}^{t}\sigma_0(X^i_s,E^i_s)dB^0_s\\
	    E^i_t &= E^i_0+\displaystyle\int_{[0,t]\times\R_+}\1_{\big\{u\leq K_{\mu_s}(X^i_s)\1_0(E^i_{s^-})+\gamma \1_1(E^i_{s^-})\big\}}\Q^i(ds,du)\\
		\mu_t &= \mathcal{L}\left((X^i_t, E^i_t)\vert\F^0_t\right).
	\end{cases}
\end{equation}
Recall that we denote by $\bar{Z}^i=(\bar{X}^i_t,\bar{E}^i_t)_{1\leq i \leq N}$ the unique solution to that system with the initial conditions $\left\{(\bar{X}^i_0,\bar{E}^i_0)\right\}_{1\leq i \leq N}$.

\subsection{Control the distance between the original system and an auxiliary system}

We can first perform a stability type estimate for the distance between the trajectories of the original system \eqref{system} and the  system of i.i.d. agents \eqref{copies}.

\begin{proof}[Proof of Theorem \ref{Prop_control_couplings}]
The first part of the proof is a slight modification of the proof for the stability of the solution to the SDE~\eqref{limit_SDE}. Indeed, we can establish an estimate similar to~\eqref{displacement_estimate} for the difference between the pairs $\bigl\{(X^{i,N},E^{i,N}), (\bar{X}^i,\bar{E}^i)\bigr\}$, $i=1,\dots,N$, in the path space $D([0,T],\Pi)$ endowed with the supremum in time distance. First, it is not hard to handle the diffusion terms imitating the proof of theorem \ref{exist_unique_thm}, the difference only arise at the two mean-field terms: the drift term and the jump term related to infection. Recall that conditionally on $\F^0$, the laws of the individuals remain the same and $\mathcal{L}\left(\big(\bar{X}^i_t, \bar{E}^i_t\big)\vert\F^0_t\right) =\mu_t$ for all $i\in\{1,\dots,N\}$.

We treat the drift terms first. Since  the  choice  of  the  coupling $m_s$ in the proof of Theorem \ref{exist_unique_thm} is  arbitrary, at the last inequality we can obtain an appropriate bound due to the definition of the 1-Wassertein metric
\begin{align*}
    \E&\Bigg[\sup_{0\leq s\leq t}\biggl|\int_0^s\frac{1}{N}\sum_{j=1}^{N}V \bigl(X^{i,N}_r,E^{i,N}_r,X^{j,N}_r,E^{j,N}_r\bigr) \,dr-\int_0^s V_{\mu_r}(\bar{X}^i_r,\bar{E}^i_r)\,dr\biggr|\Bigg]\\
    &\leq\E \left[\int_{0}^{t}\left|V_{\mu^N_s}(X^{i,N}_s,E^{i,N}_s)-V_{\mu_s}(\bar{X}^i_s,\bar{E}^i_s)\right|ds\right]\\
	&\leq C\int_{0}^{t}\E\Big[|X^{i,N}_s-\bar{X}^i_s|+|E^{i,N}_s-\bar{E}^i_s|+W_1\left(\mu^N_s,\mu_s\right)\Big] ds.
\end{align*}
For the mean-field term related to the jump, we use the same coupling and apply the same argument. To be more precise, one has

\begin{align*}
	\E&\left[\underset{0\leq s\leq t}{\sup}\int_{[0,s]\times\R_+}\left|\1_{\big\{ u \le K^{i,N}_{r} \1_0(E^{i,N}_{r^-}) + \gamma \1_1(E^{i,N}_{r^-})\big\}}-\1_{\big\{u\leq K_{\mu_r}(\bar{X}^i_{r^-})\1_0(\bar{E}^i_{r^-})+\gamma \1_1(\bar{E}^i_{r^-})\big\}}\right|\Q(dr,du)\right]\\
	&\leq\E\left[\int_0^t\left|\left(K^{i,N}_s \1_0(E^{i,N}_s) + \gamma \1_1(E^{i,N}_s)\right)-\left(K_{\mu_s}(\bar{X}^i_s)\1_0(\bar{E}^i_s)+\gamma \1_1(\bar{E}^i_s)\right)\right|ds\right]\\
	&\leq \E\left[\int_0^t\left| K^{i,N}_s \1_0(E^{i,N}_s)-K_{\mu_r}(\bar{X}^i_s)\1_0(\bar{E}^i_s)\right|ds\right]+\gamma\E\left[\int_0^t\left|\1_1(E^{i,N}_s)-\1_1(\bar{E}^i_s)\right|ds\right]\\
	&=I_1+I_2
\end{align*}
We will treat $I_1$ which contains the mean-field term $K^{i,N}_{s}$,

\begin{align*}
			I_1 =&\E\left[\int_0^t\left|\int_{\Pi}K(X^{i,N}_s,x)\1_0(E^{i,N}_s)\1_1(e)\mu^N_s(dx,de)-\int_{\Pi}K(\bar{X}^i_s,y)\1_0(\bar{E}^i_s)\1_1(f)\mu_s(dy,df)\right|ds\right]\\
% 			\leq&\E\left[\int_0^t\left|\int_{\Pi\times\Pi}\Bigl(K(X^{i,N}_s,x)-K(\bar{X}^i_s,y)\Bigr)\1_0(E^{i,N}_s)\1_1(e)m_s(dx,de,dy,df)\right|ds\right]\\
% 			&+\E\left[\int_0^t\left|\int_{\Pi\times\Pi}K(\bar{X}^i_s,y)\Bigl(\1_0(E^{i,N}_s)\1_1(e)-\1_0(\bar{E}^i_s)\1_1(f)\Bigr)m_s(dx,de,dy,df)\right|ds\right]\\
% 			\leq& L_K\int_{0}^{t}\E\left[|X^{i,N}_s-\bar{X}^i_s|+\int_{\Pi\times\Pi}|x-y|m_s(dx,de,dy,df)\right]ds\\
% 			&+\|K\|_{\infty}\int_0^t\E\left[|E^{i,N}_s-\bar{E}^i_s|+\int_{\Pi\times\Pi}|e-f|m_s(dx,de,dy,df)\right]ds\\
			\leq& C\int_{0}^{t}\E\Big[|X^{i,N}_s-\bar{X}^i_s|+|E^{i,N}_s-\bar{E}^i_s|+W_1\left(\mu^N_s,\mu_s\right)\Big] ds.
		\end{align*}
The remaining terms are treated similarly as in the proof of Theorem \ref{exist_unique_thm}. Combining them we can deduce the following bound,
\begin{align*}
\begin{split}
    \E &\left[ \sup_{0 \leq s \leq t} \left|X^{i,N}_s-\bar{X}^i_s\right|+|E^{i,N}_s-\bar{E}^i_s|\right]\\
    &\leq 2\E \left[|X^{i,N}_0-\bar{X}^i_0|+|E^{i,N}_0-\bar{E}^i_0|\right]+C\int_{0}^{t}\E\Big[|X^{i,N}_s-\bar{X}^i_s|+|E^{i,N}_s-\bar{E}^i_s|+W_1\left(\mu^N_s,\mu_s\right)\Big]ds.
\end{split}
\end{align*}
On the other hand, one has
\begin{equation}\label{estimate_for_couplings_1}
    \begin{aligned}
	\E\left[W_1\left(\mu^N_s,\mu_s\right)\right]=&\E\left[W_1\left(\frac{1}{N}\sum_{i=1}^{N}\delta_{Z^{i,N}_s},\mu_s\right)\right]\\
	\leq& \E\left[W_1\left(\frac{1}{N}\sum_{i=1}^{N}\delta_{Z^{i,N}_s},\frac{1}{N}\sum_{i=1}^{N}\delta_{\bar{Z}^i_s}\right)\right]+\E\left[W_1\left(\frac{1}{N}\sum_{i=1}^{N}\delta_{\bar{Z}^i_s},\mu_s\right)\right]\\
	\leq& \frac{1}{N}\sum_{i=1}^{N}\E\left[|X^{i,N}_s-\bar{X}^i_s|+|E^{i,N}_s-\bar{E}^i_s|\right]+\E\left[W_1\left(\frac{1}{N}\sum_{i=1}^{N}\delta_{\bar{Z}^i_s},\mu_s\right)\right].
	\end{aligned}
\end{equation}
From this and the fact that $(Z^{i,N},\bar{Z}^i), i=1,\dots,N$ are identically distributed, we have the following
\begin{align*}
    \E\bigg[\sup_{0 \leq s \leq t}& \left|X^{i,N}_s-\bar{X}^i_s\right|+|E^{i,N}_s-\bar{E}^i_s|\bigg]\\
    \leq& 2\E \left[|X^{i,N}_0-\bar{X}^i_0|+|E^{i,N}_0-\bar{E}^i_0|\right]\\
    &+C\int_{0}^{t}\E\left[|X^{i,N}_s-\bar{X}^i_s|+|E^{i,N}_s-\bar{E}^i_s|+W_1\left(\frac{1}{N}\sum_{i=1}^{N}\delta_{\bar{Z}^i_s},\mu_s\right)\right]ds.
\end{align*}
Now applying the Gronwall lemma, we deduce that
\begin{align*}
    \begin{split}
        \E&\left[ \sup_{0 \leq s \leq t} \left|X^{i,N}_s-\bar{X}^i_s\right|+|E^{i,N}_s-\bar{E}^i_s|\right]\\
    &\leq e^{Ct}\left(2\E\left[|X^{i,N}_0-\bar{X}^i_0|+|E^{i,N}_0-\bar{E}^i_0|\right]+\int_{0}^{t}\E\left[W_1\left(\frac{1}{N}\sum_{i=1}^{N}\delta_{\bar{Z}^i_s},\mu_s\right)\right]ds\right).
    \end{split}
\end{align*}
Since conditionally upon $\F^0$, $\{\bar{Z}^i_t=(\bar{X}^i_t,\bar{E}^i_t)\}_{i=1,\dots,N}$ are $N$ i.i.d.r.v., by the Glivenko-Cantelli theorem, their empirical measure converges $a.s$ to $\mu_t$. So the second term in the r.h.s. in the above bound goes to $0$. Moreover, if $\E\left[|Z^i_t|^q\right]<\infty$ for some $q>2$, it is possible to give a rate of convergence measured in the Wasserstein distance. Indeed, we are able to adapt to the conditional case the crucial result obtained by Fournier-Guillin~\cite{Fournier-Guillin}. See Proposition~\ref{app:FG_new} in the Appendix for a precise statement.

% 	\begin{equation*}
% 	\E\left[W_1\left(\frac{1}{N}\sum_{j=1}^{N}\delta_{\bar{Z}^j_t},\mu_t\right)\right]\leq C\E\left[|Z^1_t|^q\right]^{1/q}
% 	\begin{cases}
% 	N^{-1/2}+N^{-(q-1)/q},& d<2
% 	\\ N^{-1/2}\log N+N^{-(q-1)/q},& d=2
% 	\\ N^{-1/d}+N^{-(q-1)/q},& d>2
% 	\end{cases}
% 	\end{equation*}
Hence by Propositions~\ref{app:FG_new} and~\ref{finite_moments}, we can deduce the desired stability estimate
\begin{equation}\label{estimate_for_couplings_2}
		\E \left[ \sup_{0 \leq s \leq t} \left|X^{i,N}_s-\bar{X}^i_s\right|+|E^{i,N}_s-\bar{E}^i_s|\right]\leq C(t)\left(\E \left[|X^{i,N}_0-\bar{X}^i_0|+|E^{i,N}_0-\bar{E}^i_0|\right]+\alpha_d(N)\right)
	\end{equation}
with $\alpha_d(N)$ as defined in the statement of Theorem~\ref{Prop_control_couplings}.
\end{proof}

As a consequence, at any time $t$, we can give a quantitative estimate of convergence of the empirical measure of the original sytem towards its expected limit $\mu_t$, measured in the Wasserstein distance.

\subsection{Proof of Corollary~\ref{rate_of_convergence}}
\begin{proof}
Combining \eqref{estimate_for_couplings_1} and \eqref{estimate_for_couplings_2}, we obtain
\begin{align*}
	\E\left[W_1\left(\mu^N_s,\mu_s\right)\right]\leq& \frac{1}{N}\sum_{i=1}^{N}\E\left[|X^{i,N}_s-\bar{X}^i_s|+|E^{i,N}_s-\bar{E}^i_s|\right]+\E\left[W_1\left(\frac{1}{N}\sum_{i=1}^{N}\delta_{\bar{Z}^i_s},\mu_s\right)\right]\\
	\leq& C(t)\left(\E\left[\frac{1}{N}\sum_{i=1}^{N}|X^{i,N}_0-\bar{X}^i_0|+|E^{i,N}_0-\bar{E}^i_0|\right]+\alpha_d(N)\right).
	\end{align*}
Thanks to the exchangeability and a coupling result, see for instance \cite[Proposition A.1]{Fournier-Hauray} or \cite[Proposition 2.14]{Hauray-Mischler}, we have the right to choose the initial conditions of the coupling such that, almost surely,

\begin{equation*}
    W_1\left(\mu^N_0,\bar{\mu}^N_0\right)=\frac{1}{N}\sum_{i=1}^{N}|X^{i,N}_0-\bar{X}^i_0|+|E^{i,N}_0-\bar{E}^i_0|,
\end{equation*}
which allows us to deduce
    \begin{align*}
		\E\left[W_1\left(\mu^N_t,\mu_t\right)\right]\leq& C(t)\left(\E\left[W_1\left(\mu^N_0,\bar{\mu}^N_0\right)\right]+\alpha_d(N)\right)\\
		\leq& C(t)\left(\E \left[W_1\left(\mu^N_0,\mu_0\right)\right]+\E\left[W_1\left(\bar{\mu}^N_0,\mu_0\right)\right]+\alpha_d(N)\right)\\
		\leq& C(t)\left( \E\left[W_1\left(\mu^N_0,\mu_0\right)\right]+\alpha_d(N)\right)
	\end{align*}
where we allow the constant to change from line to line, and use the quantitative bound for $W_1\left(\bar{\mu}^N_0,\mu_0\right)$ for the last inequality.
\end{proof}

\section{Appendix}

\subsection{Proof of the equivalences in Definition \ref{def:cond-chaos}} \label{app:cond_chaos}

\begin{proof}
$(i) \Rightarrow (ii)$ is obvious.

$(ii)\Rightarrow (iii)$. Take any $\phi \in C_b(\Pi)$,
\begin{align*}
\E\left[\left<\mu^N-\mu, \phi\right>^2 \, \bigg\vert \, \F^0\right]=&\dfrac{1}{N^2}\sum_{i,j=1}^N\E\left[\phi(Z^{i,N})\phi(Z^{j,N})\vert \F^0\right] \\
&-\dfrac{2}{N}\sum_{i=1}^N\E\left[\phi(Z^{i,N})\left<\mu,\phi\right>\vert \F^0\right]+
%\left<\mu,\phi\right>^2
\E\bigl[\left<\mu,\phi\right>^2\,\big\vert \, \F^0\bigr]
\end{align*}
Using the conditional exchangeability we rewrite the r.h.s. as
\begin{align*}
\dfrac{1}{N}\E\left[\phi(Z^{1,N})^2 \,\vert \, \F^0\right]+\dfrac{N-1}{N}\E\left[\phi(Z^{1,N})\phi(Z^{2,N}) \,\vert \,\F^0\right]
-2 \left<\mu,\phi\right> \E\left[\phi(Z^{1,N})   \, \Big\vert \, \F^0\right]+
\left<\mu,\phi\right>^2
%\E\left[\left<\mu,\phi\right>^2\ \, \big\vert \, \F^0\right]
\end{align*}
which tends to $0$ by (ii). Therefore,
		\[\E\left[\left<\mu^N,\phi\right> \, \vert \, \F^0\right]\to\left<\mu,\phi\right>
		\]
Since $\mu$ is $\F^0$-measurable and the above statement is true for any bounded continuous $\phi$, we conclude that $\mu^N$ converges weakly to $\mu$ $\Pb^0$-\text{a.s.}

\medskip
$(iii) \Rightarrow (i)$.
We consider only the case $k=2$ of that implication, but the general case $k \ge2$ could be handled in a similar way.

Let $\phi_1,\phi_2 \in C_b(\Pi)$. We will prove the convergence in law of $(Z^{1,N},Z^{2,N})$ using only functions of the form $\phi(z_1,z_2) = \phi_1(z_1) \phi_2(z_2)$, whose linear combinations are dense in $C_p(\Pi^2)$.  By the triangular inequality, one has
\begin{align}
	\bigg|\E\Big[\phi_1&(Z^{1,N})\phi_2(Z^{2,N}) \, \Big\vert \,\F^0\Big]-\left<\mu,\phi_1\right>\left<\mu,\phi_2\right>\bigg| \nonumber \\
	\leq&\bigg|\E\Big[\phi_1(Z^{1,N})\phi_2(Z^{2,N}) \, \Big\vert \,\F^0\Big]-\E\Big[\left<\mu^N,\phi_1\right>\left<\mu^N,\phi_2\right> \, \Big\vert \, \F^0\Big]\bigg|\label{iii_ii_bound_1}\\
	&+\bigg|\E\Big[\left<\mu^N,\phi_1\right>\left<\mu^N,\phi_2\right>
	 \, \Big\vert \,\F^0\Big]-\E\Big[\left<\mu,\phi_1\right>\left<\mu,\phi_2\right> \, \Big\vert\, \F^0\Big]\bigg|  \label{iii_ii_bound_2}
\end{align}
Using the conditional exchangeability, we have an upper bound for \eqref{iii_ii_bound_1},
\begin{align*}
	\bigg|\E\Big[\phi_1&(Z^{1,N})\phi_2(Z^{2,N})\vert\F^0\Big]-\E\Big[\left<\mu^N,\phi_1\right>\left<\mu^N,\phi_2\right>\vert\F^0\Big]\bigg|
	\\
	=&\Bigg|\frac{1}{N(N-1)}\sum_{i,j=1,i\neq j}^{N}\E\left[\phi_1(Z^{i,N})\phi_2(Z^{j,N}) \vert \F^0\right]
	-\frac{1}{N^2}\sum_{i,j=1}^{N}\E\left[\phi_1(Z^{i,N})\phi_2(Z^{j,N}) \vert \F^0\right]\Bigg|\\
	\leq&\Bigg|\left(\frac{1}{N(N-1)}-\frac{1}{N^2}\right)\sum_{i,j=1,i\neq j}^{N}\E\Big[\phi_1(Z^{i,N})\phi_2(Z^{j,N})\vert\F^0\Big]\Bigg|\\
	& \hspace{0.5cm} +\Bigg|\frac{1}{N^2}\sum_{i=1}^{N}\E\left[\phi_1(Z^{i,N})\phi_2(Z^{i,N})\vert \F^0\right]\Bigg|\\
	\leq& \frac1N \|\phi_1\|_{\infty}\|\phi_2\|_{\infty}+\frac1N \|\phi_1\|_{\infty}\|\phi_2\|_{\infty}
\end{align*}
where the last quantity tends to $0$ as $N\rightarrow \infty$.

Considering \eqref{iii_ii_bound_2}, we can prove that it converges to $0$, using the point $(iii)$ and the fact  that the function 
$\tilde \mu \mapsto \langle \ \tilde \mu, \phi_1\rangle \langle \tilde \mu, \phi_2 \rangle$
is continuous on $\mathcal P(\Pi)$.

\iffalse
 in $L^1$,
\begin{align*}
	\E\Big[ & \left|\E[\left<\mu^N,\phi_1\right>\left<\mu^N,\phi_2\right>\vert\F^0]-\E[\left<\mu,\phi_1\right>\left<\mu,\phi_2\right>\vert\F^0]\right|\Big]\\
	&\leq\E\Big[\left|\left<\mu^N,\phi_1\right>-\left<\mu,\phi_1\right>\right|\left|\left<\mu^N,\phi_2\right>\right|\Big]+\E\Big[\left|\left<\mu^N,\phi_2\right>-\left<\mu,\phi_2\right>\right|\left|\left<\mu,\phi_1\right>\right|\Big]\\
	&\leq \|\phi_2\|_{\infty}\E\Big[\left|\left<\mu^N,\phi_1\right>-\left<\mu,\phi_1\right>\right|\Big]+\|\phi_1\|_{\infty}\E\Big[\left|\left<\mu^N,\phi_2\right>-\left<\mu,\phi_2\right>\right|\Big]\rightarrow 0 \qquad as \; N\rightarrow\infty.
\end{align*}
the last line follows from $(iii)$.
\fi
\end{proof}

\subsection{Proof of Lemma \ref{lemma:cond-exchangeability}}
\label{app:cond_exchangeability}

\begin{proof}
The sytem of  equation~\eqref{system} is symetric under any permutation of the individuals and their individual noises. Hence, if $Z^N_t:=\left(Z^{1}_t,\dots, Z^N_t\right)$ is a strong solution of the system with initial condition $\left(Z^1_0, \dots, Z^N_0\right)$ then $Z^{\sigma, N}_t:=\left(Z_{t}^{\sigma(1), N}, \dots, Z_{t}^{\sigma(N), N}\right)$ is also a strong solution of the equation with initial condition $\left(Z^{\sigma(1)}_0, \dots, Z^{\sigma(N)}_0\right)$ with permuted individual noise. Moreover, by the assumption that the initial laws of $Z^{N}_0$ and $Z^{\sigma, N}_0$ are equal, the conclusion follows using the weak uniqueness of solutions to the SDE system~\eqref{system} with jumps. In fact, the standard approach of \cite[Theorem 9.1]{Ikeda-Watanabe} applies to our particular case, since all the coefficients appearing in the system~\eqref{system} are regular enough.

% More precisely for all $t\geq 0$, one has
% \[
% \bigl(Z_{t}^{N},B_{t}^{0}\bigr)_{\#} \mathbb{P}=\bigl(Z_{t}^{\pi, N},B_{t}^{0}\bigr)_{\#} \mathbb{P}
% \]
% Therefore, for all $A \in \mathcal{F}_{t}^{0}$ and for every $\phi \in C_{b}\left(\left(\mathbb{R}^{d}\right)^{N}\right)$,
% \[
% \mathbb{E}\left[\1_{A} \phi\left(Z_{t}^{N}\right)\right] = \mathbb{E}\left[\1_{A} \phi\left(Z_{t}^{\pi, N}\right)\right]
% \]
% Because the expectation of $\phi\left(Z_{t}^{N}\right)$ and $\phi\left(Z_{t}^{\pi, N}\right)$ coincide on every element of a basis of $\mathcal{F}_{t}^{0}$, their conditional expectation coincide as well.

\end{proof}

\subsection{Rate of convergence in the empirical law of large number: the conditionnally independence case} \label{app:FG_new}

\begin{proposition} \label{prop:FG_new}
For any $d \in \N$ and $q \in (2, \infty)$, there exists a constant $C_{d,q}$ such that the following holds. For any $\mathbb P_0$-measurable random probability $\mu$ on $\Pi$ satisfying
\[
\E \biggl[ \int_{\R^d} |z|^q \mu(dz) \biggr] 
= \E \bigl[|Z_i|^q \bigr] 
< \infty ,
\]
if conditionnally upon $\F_0$, the $\Pi$-valued random variables $(Z_i)_{i \ge 1}$ are i.i.d. with conditional law $\mu$, and denoting the random empirical measure by $\mu^N_Z= \frac{1}{N}\sum_{i=1}^{N}\delta_{Z_i}$, we have
\begin{equation*}
	\E\left[W_1\left(\mu^N_Z ,\mu \right)\right]\leq C_{d,q}
	\E\left[|Z_1|^q\right]^{1/q}
	\begin{cases}
	N^{-1/2},& d =1 %+N^{-(q-1)/q} 
	\\ N^{-1/2}\log N,& d=2 % +N^{-(q-1)/q}
	\\ N^{-1/d},& d \ge 3 % +N^{-(q-1)/q}
	\end{cases}
\end{equation*}
\end{proposition}

\begin{proof}
For a given realisation of the common noise, we rely first on a result by Fournier Guillin~\cite[Theorem 1]{Fournier-Guillin} (which is the best result so far after a long sequence of previous partial or less  accurate result) in the case of i.i.d. random variables with a given (deterministic) common law $\mu$.
Applying it (with $p=1$) in our special case\footnote{remark that the particular assumption on $q$ in\cite[Theorem 1]{Fournier-Guillin} is fulfilled by our particular choice $q>2$.}, it gives the following estimate conditionnally on $\F^0$,
\[\Pb^0-\text{a.s.}, \quad
\E\left[W_1\left(\mu^N_Z ,\mu \right) \,\big\vert\, \F^0 \right] \leq 
C_{d,q}
	\E\left[|Z_1|^q \,\big\vert\, \F^0 \right]^{1/q}
	\begin{cases}
	N^{-1/2}+N^{-(q-1)/q},& d =1
	\\ N^{-1/2}\log N+N^{-(q-1)/q},& d=2
	\\ N^{-1/d}+N^{-(q-1)/q},& d \ge 3
	\end{cases}
\]
First, we could remove all the second power of $N$ involving $q$ in the three cases since it is always smaller that the first term when $q>2$.

Taking then the expectation w.r.t. to $\mathbb P_0$ in the above inequality, and using since $q>1$ the Jensen inequality on the $q$ moment
\[
\E \Bigl[\E\left[|Z_1|^q \,\big\vert\, \F^0 \right]^{1/q}
\Bigr]\leq \E\left[|Z_1|^q \right]^{1/q}
\]
we get the desired result.
\end{proof}

\end{document}